\documentclass[11pt,english,journal,onecolumn,12pt]{IEEEtran}
\usepackage[T1]{fontenc}
\usepackage[latin9]{inputenc}
\usepackage[letterpaper]{geometry}
\usepackage{amsthm}
\usepackage{amsmath}
\usepackage{amssymb}
\usepackage{graphicx}
\usepackage{setspace}
\makeatletter
\theoremstyle{plain}
\newtheorem{thm}{\protect\theoremname}
\theoremstyle{plain}
\newtheorem{prop}[thm]{\protect\propositionname}
\theoremstyle{plain}
\newtheorem{lem}[thm]{\protect\lemmaname}
\theoremstyle{plain}
\newtheorem{cor}[thm]{\protect\corollaryname}

\usepackage{amsfonts}%\usepackage{dropping}

\usepackage{babel}

\makeatother

\usepackage{babel}
\providecommand{\corollaryname}{Corollary}
\providecommand{\lemmaname}{Lemma}
\providecommand{\propositionname}{Proposition}
\providecommand{\theoremname}{Theorem}

\begin{document}

\title{On the Optimal Scheduling of Independent, Symmetric and Time-Sensitive
Tasks}

\author{Fabio Iannello$^{1}$, Osvaldo Simeone$^{2}$ and Umberto Spagnolini$^{3}$%
\thanks{$^{1,3}$F. Iannello (corresponding author: iannello@elet.polimi.it,
phone: +390223993604, fax +390223993413) and U. Spagnolini (spagnoli@elet.polimi.it)
are with Politecnico di Milano, P.zza L. da Vinci 32, 20133 Milan,
Italy.$^{1,2}$F. Iannello and O. Simeone (osvaldo.simeone@njit.edu)
are with the CWCSPR, New Jersey Institute of Technology, U. Heights,
Newark, NJ, 07102, USA %
}}
\maketitle
\begin{abstract}
Consider a discrete-time system in which a \emph{centralized controller}
(CC) is tasked with assigning at each time interval (or \emph{slot})
$K$ resources (or \emph{servers}) to $K$ out of $M\geq K$\emph{
nodes}. When assigned a server, a node can execute a \emph{task. }The
tasks are independently generated at each node by stochastically symmetric
and memoryless random processes and stored in a finite-capacity\emph{
task queue}. Moreover, they are \emph{time-sensitive} in the sense
that within each slot there is a non-zero probability that a task
expires before being scheduled. The scheduling problem is tackled
with the aim of maximizing the number of tasks completed over time
(or the \emph{task-throughput}) under the assumption that the CC has
no direct access to the state of the task queues. The scheduling decisions
at the CC are based on the outcomes of previous scheduling commands,
and on the known statistical properties of the task generation and
expiration processes.

Based on a Markovian modeling of the task generation and expiration
processes, the CC scheduling problem is formulated as a partially
observable Markov decision process (POMDP) that can be cast into the
framework of restless multi-armed bandit (RMAB) problems. When the
task queues are of capacity one, the optimality of a myopic (or greedy)
policy is proved. It is also demonstrated that the MP coincides with
the Whittle index policy. For task queues of arbitrary capacity instead,
the myopic policy is generally suboptimal, and its performance is
compared with an upper bound obtained through a relaxation of the
original problem.

Overall, the settings in this paper provide a rare example where a
RMAB problem can be explicitly solved, and in which the Whittle index
policy is proved to be optimal. 
\end{abstract}

\section{Introduction and System Model\label{sec:Introduction}}

\begin{figure}[h!]
\centering \includegraphics[width=4.5in]{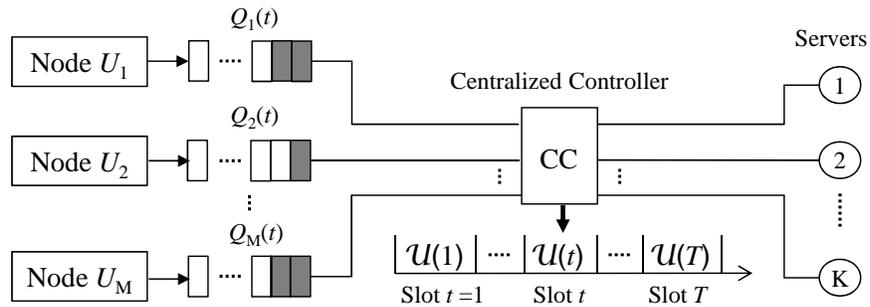}
\caption{The centralized controller (CC) assigns $K$ resources (servers) to
$K$ out of $M\geq K$ nodes to complete their tasks in each slot
$t$. The tasks of node $U_{i}$ at slot $t$ are stored in a task
queue $Q_{i}(t).$}

\label{fig:System_model}
\end{figure}

The problem of scheduling concurrent tasks under resource constraints
finds applications in a variety of fields including communication
networks \cite{art:Data_networks}, distributed computing \cite{art:Grid_computing_2}
and virtual machine scenarios \cite{art:Virtual_Machine}. In this
paper we consider a specific instance of this general problem in which
a\emph{ centralized controller} (CC) is tasked with assigning at each
time interval (or \emph{slot}) $K$ resources, referred to as \emph{servers},
to $K$ out of $M\geq K$ \emph{nodes} as shown in Fig. \ref{fig:System_model}.
A server can complete a single task per\emph{ }slot and can be assigned
to one node per time interval. The tasks are generated at the $M$
nodes by stochastically symmetric, independent and memoryless random
processes. The tasks are stored by each node in a finite-capacity
\emph{task queue}, and they are \emph{time-sensitive} in the sense
that at each slot there is a non-zero probability that a task expires
before being completed successfully. It is assumed that the CC has
no direct access to the node queues, and thus it is not fully informed
of their actual states. Instead, the scheduling decision is based
on the outcomes of previous scheduling commands, and on the statistical
knowledge of the task generation and expiration processes. If a server
is assigned to a node with an empty queue, it remains idle for the
whole slot. The purpose here is thus to pair servers to nodes so as
to maximize the average number of successfully completed tasks within
either a finite or infinite number of slots (\emph{horizon}), which
we refer to as \emph{task-throughput}, or simply throughput. 

\subsection{Markov Formulation\label{sub:Markov-Formulation}}

\begin{figure}[h!]
\centering \includegraphics[clip,width=12cm]{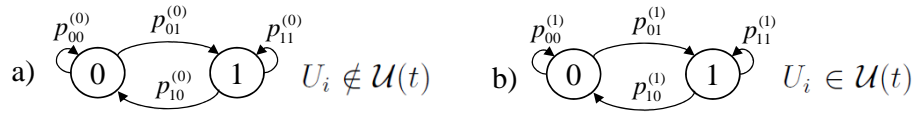}
\caption{Markov model for the evolution of the state of the task queue $Q_{i}(t)\in\{0,1\}$,
when the node $U_{i}$: a) is not scheduled in slot $t$ (i.e., $U_{i}\notin\mathcal{U}(t)$);
b) is scheduled in slot $t$ (i.e., $U_{i}\in\mathcal{U}(t)$).}

\label{fig:Full_DMC}
\end{figure}

We now introduce the stochastic model that describes the evolution
of the task queues across slots. In this section we consider task
queues of capacity one (see Sec. \ref{sec:Extension_capacity} for
capacity larger than one), where $Q_{i}(t)\in\{0,1\}$ denotes the
number of tasks in the queue of node $U_{i}$, for $i\in\{1,...,M\}$.
The stochastic evolution of queue $Q_{i}(t)$ is shown in Fig. \ref{fig:Full_DMC}
as a function of the scheduling decision $\mathcal{U}(t)$, which
consists in the assignment at each slot $t$ of the $K$ servers to
a subset $\mathcal{U}(t)\subseteq\{U_{1},...,U_{M}\}$ of $K$ nodes,
with $\left|\mathcal{U}(t)\right|=K$. 

At each slot, node $U_{i}$ can be either scheduled ($U_{i}\in\mathcal{U}(t)$)
or not ($U_{i}\notin\mathcal{U}(t)$). If $U_{i}$ is not scheduled
(i.e., $U_{i}\notin\mathcal{U}(t)$, see Fig. \ref{fig:Full_DMC}-a))
and there is a task in its queue (i.e., $Q_{i}(t)=1$), then the task
expires with probability (w.p.) $p_{10}^{(0)}=\Pr[Q_{i}(t+1)=0|Q_{i}(t)=1,U_{i}\notin\mathcal{U}(t)]$,
while it remains in the queue w.p. $p_{11}^{(0)}=1-p_{10}^{(0)}$.
Instead, if node $U_{i}$ is scheduled (i.e., $U_{i}\in\mathcal{U}(t)$,
see Fig. \ref{fig:Full_DMC}-b)) and $Q_{i}(t)=1$, its task is completed
successfully and its queue in the next slot is either empty or full
w.p. $p_{10}^{(1)}=\Pr[Q_{i}(t+1)=0|Q_{i}(t)=1,U_{i}\in\mathcal{U}(t)]$
and $p_{11}^{(1)}=1-p_{10}^{(1)}$, respectively. Probability $p_{11}^{(1)}$
accounts for the possible arrival of a new task. If $Q_{i}(t)=0$
the probabilities of receiving a new task when $U_{i}$ is not scheduled
and scheduled are $p_{01}^{(0)}=\Pr[Q_{i}(t+1)=1|Q_{i}(t)=0,U_{i}\notin\mathcal{U}(t)]$
and $p_{01}^{(1)}=\Pr[Q_{i}(t+1)=1|Q_{i}(t)=0,U_{i}\in\mathcal{U}(t)]$,
respectively, while the probabilities of receiving no task are $p_{00}^{(0)}=1-p_{01}^{(0)}$
and $p_{00}^{(1)}=1-p_{01}^{(1)}$, respectively.

\subsection{Related Work and Contributions\label{sub:Related-Work-and}}

In this work we assume that the CC has no direct access to the state
of the task queues $Q_{1}(t),...,Q_{M}(t)$, while it knows the transitions
probabilities $p_{xy}^{(u)}$, with $x,y,u\in\{0,1\}$, and the outcomes
of previously scheduled tasks. The scheduling problem is thus formalized
as a partially observable Markov decision process (POMDP) \cite{art:Monahan},
and then cast into a restless multi-armed bandit (RMAB) problem \cite{book:Gittins}.
A RMAB is constituted by a set of arms (the queues in our model),
a subset of which needs to be activated (or scheduled) in each slot
by the controller. 

To elaborate, we assume that the transition probabilities of the Markov
chains in Fig. \ref{fig:Full_DMC}, the number of nodes $M$ and servers
$K$ are such that\begin{subequations}
\begin{align}
 & m=M/K,\text{ }\text{is integer, and}\label{eqn: integer_assumptions}\\
 & p_{11}^{(1)}\leq p_{01}^{(1)}\leq p_{01}^{(0)}\leq p_{11}^{(0)}.\label{eqn: inequality_transition_probabilities}
\end{align}
\label{eq:RR_conditions}\end{subequations}Assumption (\ref{eqn: integer_assumptions})
states that the ratio $m=M/K$ between the numbers $M$ of nodes and
$K$ of servers is an integer, generalizing the single-server case
($K=1$). Proving the results provided later in this paper for the
case of non-integer $m$ remains an open problem. Assumption (\ref{eqn: inequality_transition_probabilities})
is motivated as follows. The inequality $p_{11}^{(1)}\leq p_{01}^{(1)}$
imposes that the probability that a new task arrives when the task
queue is full and the node is scheduled ($p_{11}^{(1)}$) is no larger
than when the task queue is empty ($p_{01}^{(1)}$). This applies,
e.g., when the arrival of a new task is independent on the queue's
state and scheduling decisions (i.e., $p_{11}^{(1)}=p_{01}^{(1)}$),
or when a new task is not accepted when the queue is full, i.e., $p_{11}^{(1)}=0$.
Inequality $p_{01}^{(1)}\leq p_{01}^{(0)}$ applies, e.g., when the
task generation process does not depend on the queue's state and on
the scheduling decisions, so that $p_{01}^{(1)}=p_{01}^{(0)}$, or
when a new task cannot be accepted while the node is scheduled even
if the queue is empty ($p_{01}^{(1)}=0$). Inequality $p_{01}^{(0)}\leq p_{11}^{(0)}$
indicates that, when a node is not scheduled, the probability $p_{01}^{(0)}$
that its task queue is full in the next slot, given that it is currently
empty, is smaller than the probability $p_{11}^{(0)}$ that the task
queue is full in the next slot given that it is currently full. This
applies, e.g., when the task generation and expiration processes are
independent of each other.

\textbf{Main Contributions:} When the task queues are of capacity
one, and under assumptions (\ref{eq:RR_conditions}), we first show
that the \emph{myopic policy} (MP) for the RMAB at hand is a round
robin (RR) strategy that: \emph{i}) re-numbers the nodes in a decreasing
order according to the initial probability that their respective task
queue is full; and then \emph{ii}) schedules the nodes periodically
in group of $K$ by exploiting the initial ordering. The MP is then
proved to be throughput-optimal. We then show that, for the special
case in which $p_{01}^{(0)}=p_{01}^{(1)}$ and $p_{10}^{(0)}=p_{11}^{(1)}=0$,
the MP coincides with the Whittle index policy, which is a generally
suboptimal index strategy for RMAB problems \cite{art:whittle}. Finally,
we extend the model of Sec. \ref{sub:Markov-Formulation} to queues
with an arbitrary capacity $C$. Characterizing optimal policies for
$C>1$ is significantly more complicated than the case of $C=1$.
Hence, inspired by the optimality of the MP for $C=1$, we compare
the performance of the MP for $C>1$, with a upper bound based on
a relaxation of the scheduling constraints of the original RMAB problem
\cite{art:whittle}. It is recalled that the results in this paper
represent a rare case in which the optimal policy for a RMAB can be
found explicitly \cite{book:Gittins}.

\textbf{Related Work:} The work in this paper is related to the works
\cite{art:Zhao-myopic,art:Ahmad-vector-myopic}, in which a RMAB problem
similar to the one in this paper is addressed. However, the main difference
between our RMAB and the one in \cite{art:Zhao-myopic,art:Ahmad-vector-myopic}
is the evolution of the arms across slots. In particular, in our RMAB,
each arm evolves across a slot depending on the scheduling decision
taken by the controller, while in \cite{art:Zhao-myopic,art:Ahmad-vector-myopic},
the evolution of the arms does not depend on the scheduling decision.
The transition probabilities for the RMAB in \cite{art:Zhao-myopic,art:Ahmad-vector-myopic}
are thus equivalent to setting $p_{01}^{(0)}=p_{01}^{(1)}$ and $p_{11}^{(0)}=p_{11}^{(1)}$
in the Markov chains of Fig. \ref{fig:Full_DMC}. For instance, our
model applies to scenarios in which the arms are, e.g., data queues,
where each arm draws a data packet from its queue only when scheduled.
Instead, the model in \cite{art:Zhao-myopic,art:Ahmad-vector-myopic}
applies to scenarios in which the arms are, e.g., communication channels,
whose quality evolves across slots regardless whether they are selected
for transmission or not.

In \cite{art:Zhao-myopic} it is shown that the MP is optimal for
$p_{01}^{(0)}=p_{01}^{(1)}\leq p_{11}^{(0)}=p_{11}^{(1)}$ with $K=1,$
while \cite{art:Ahmad-vector-myopic} extends this result to an arbitrary
$K.$ The work \cite{art:Zhao-myopic} also demonstrates that the
MP in not generally optimal in the case $p_{01}^{(0)}=p_{01}^{(1)}\geq p_{11}^{(0)}=p_{11}^{(1)}$.
Finally, paper \cite{art:Zhao-whittle} proves the optimality of the
Whittle index policy for $p_{11}^{(0)}=p_{11}^{(1)}\leq p_{11}^{(0)}=p_{11}^{(1)}.$
We emphasize that neither our model nor the one considered in \cite{art:Zhao-myopic,art:Ahmad-vector-myopic}
subsumes the other, and the results here and in the mentioned previous
works should be considered as complementary.

\textit{Notation}: Vectors are denoted in bold, while the corresponding
non-bold letters denote the vectors components. Given a vector $\mathbf{x=[}x_{1},...,x_{M}\mathbf{]}$
and a set $\mathcal{S=}\left\{ i_{1},...,i_{K}\right\} \subseteq\left\{ 1,...,M\right\} $
of cardinality $K\leq M\mathbf{,}$ we define vector $\mathbf{x}_{\mathcal{S}}=[x_{i_{1}},...,x_{i_{K}}]$,
where $i_{1}\leq...\leq i_{K}$. A function $f(\mathbf{x})$ of vector
$\mathbf{x}$ is also denoted as $f(x_{1},...,x_{M})$ or as $f(x_{1},...,x_{l},\mathbf{x}_{\{l+1,...,M\}})$
for some $1\leq l\leq M$, or similar notations depending on the context.
Given a set $\mathcal{A}$ and a subset $\mathcal{B\subseteq A},$
$\mathcal{B}^{c}$ represents the complement of $\mathcal{B}$ in
$\mathcal{A}$.

\section{Problem Formulation\label{sec:System-Model}}

Here we formalize the scheduling problem of Sec. \ref{sec:Introduction}
(see Fig. \ref{fig:System_model}), in which the task generation and
expiration processes are modeled, independently at each node, by the
Markov models of Sec. \ref{sub:Markov-Formulation} with queues of
capacity one. Extension to task queues of arbitrary capacity is addressed
in Sec. \ref{sec:Extension_capacity}.

\subsection{Problem Definition\label{sub:Problem-Definition}}

The scheduling problem at the CC is addressed in a finite-horizon
scenario over slots $t\in\{1,...,T\}$. Let $\mathbf{Q}(t)=\left[Q_{1}(t),...,Q_{M}(t)\right]$
be the vector collecting the states of the task queue at slot $t$.
At slot $t=1,$ the CC is only aware of the initial probability distribution
$\boldsymbol{\omega}(1)=[\omega_{1}(1),...,\omega_{M}(1)]$ of $\mathbf{Q}(1)$,
whose $i$th entry is $\omega_{i}(1)=\Pr[Q_{i}(1)=1]$. Thus, the
subset $\mathcal{U}(1)$ of $\left|\mathcal{U}(1)\right|=K$ nodes
scheduled at slot $t=1$ is chosen as a function of the initial distribution
$\boldsymbol{\omega}(1)$ only. For any node $U_{i}\in\mathcal{U}(t)$
scheduled at slot $t$, an \emph{observation} is made available to
the CC at the end of the slot, while no observations are available
for non-scheduled nodes $U_{i}\notin\mathcal{U}(t)$. Specifically,
if $Q_{i}(t)=1$ and $U_{i}\in\mathcal{U}(t)$, the task of $U_{i}$
is served within slot $t$, and the CC observes that $Q_{i}(t)=1$.
Conversely, if $Q_{i}(t)=0$ and $U_{i}\in\mathcal{U}(t)$, no tasks
are completed and the CC observes that $Q_{i}(t)=0$. We define $\mathcal{O}(t)=\left\{ Q_{i}(t):\: U_{i}\in\mathcal{U}(t)\right\} $
as the set of (new) observations available at the CC at the end of
slot $t$. At time $t,$ the CC hence knows the history of all decisions
and previous observations and the initial distribution $\boldsymbol{\omega}(1)$,
namely $\mathcal{H}(t)=\left\{ \mathcal{U}(1),...,\mathcal{U}(t-1),\mathcal{O}(1),...,\mathcal{O}(t-1),\boldsymbol{\omega}(1)\right\} $,
with $\mathcal{H}(1)=\left\{ \boldsymbol{\omega}(1)\right\} $.

Since the CC has only partial information about the system state $\mathbf{Q}(t)$,
through $\mathcal{O}(t)$, the scheduling problem at hand can be modeled
as a POMDP. It is well-known that a sufficient statistics for taking
decisions in such POMDP is given by the probability distribution of
$\mathbf{Q}(t)$ conditioned on the history $\mathcal{H}(t)$ \cite{art:Cassandra},
referred to as \textit{belief}, and represented by the vector $\boldsymbol{\omega}(t)=[\omega_{1}(t),...,\omega_{M}(t)]$,
with $i$th entry given by
\begin{equation}
\omega_{i}(t)=\Pr\left[Q_{i}(t)=1|\mathcal{H}(t)\right].\label{eqn: belief_vector_entries}
\end{equation}
Since the belief $\boldsymbol{\omega}(t)$ fully summarizes the entire
history $\mathcal{H}(t)$ of past actions and observations \cite{art:Cassandra},
a \emph{scheduling policy} $\pi\mathbf{=}\left[\mathcal{U}^{\pi}(1),...,\mathcal{U}^{\pi}(T)\right]$
is defined as a collection of functions $\mathcal{U}^{\pi}(t)$ that
map the belief $\boldsymbol{\omega}(t)$ to a subset $\mathcal{U}(t)$
of $\left|\mathcal{U}(t)\right|=K$ nodes, i.e., $\mathcal{U}^{\pi}(t)$:
$\boldsymbol{\omega}(t)\rightarrow\mathcal{U}(t)$. We will refer
to $\mathcal{U}^{\pi}(t)$ as the subset of scheduled nodes, even
though, strictly speaking, it is the mapping function defined above.
The transition probabilities over the belief space are derived in
Sec. \ref{Sec: Transition_probabilities}.

The \emph{immediate reward} $R(\boldsymbol{\omega},\mathcal{U})$,
accrued by the CC when the belief vector is $\boldsymbol{\omega}$
and action $\mathcal{U}$ is taken, measures the average number of
tasks completed within the current slot, and it is
\begin{equation}
R(\mathbf{\boldsymbol{\omega}},\mathcal{U})=\sum_{U_{i}\in\mathcal{U}}\omega_{i}.\label{eqn: reward_belief_MDP_line_1}
\end{equation}
Notice that $R(\mathbf{\boldsymbol{\omega}},\mathcal{U})\leq K$ since
there are only $K$ servers.

The \emph{throughput} measures the average number of tasks completed
over the slots $\{1,...,T\}$ that, by exploiting (\ref{eqn: reward_belief_MDP_line_1})
and under policy $\pi$, is given by
\begin{equation}
V_{1}^{\pi}\left(\mathbf{\boldsymbol{\omega}}(1)\right)=\sum_{t=1}^{T}\beta^{t-1}\mathrm{E}^{\mathbf{\pi}}\left[R\left(\mathbf{\boldsymbol{\omega}}(t)\mathbf{,}\mathcal{U}^{\mathbf{\pi}}(t)\right)|\mathbf{\boldsymbol{\omega}}(1)\right].\label{eqn: discounted_avg_reward_with_BELIEF}
\end{equation}
In (\ref{eqn: discounted_avg_reward_with_BELIEF}), the expectation
$\mathrm{E}^{\mathbf{\pi}}[\cdot|\mathbf{\boldsymbol{\omega}}(1)]$,
under policy $\pi$ for initial belief $\mathbf{\boldsymbol{\omega}}(1)$,
is intended with respect to the distribution of the Markov process
$\mathbf{\boldsymbol{\omega}}(t)$, as obtained from the transition
probabilities to be derived in Sec. \ref{Sec: Transition_probabilities}.
For generality, the definition (\ref{eqn: discounted_avg_reward_with_BELIEF})
includes a discount factor $0\leq\beta\leq1$ \cite{art:Zhao-myopic},
while the infinite horizon scenario (i.e., $T\rightarrow\infty$)
will be discussed in Sec. \ref{Sec: infinite_horizon}.

The goal is to find a policy $\pi^{\ast}=\left[\mathcal{U}^{*}(1),...,\mathcal{U}^{*}(T)\right]$
that maximizes the throughput (\ref{eqn: discounted_avg_reward_with_BELIEF})
so that
\begin{eqnarray}
V_{1}^{\mathbf{\ast}}\left(\boldsymbol{\omega}(1)\right)=V_{1}^{\pi^{\ast}}\left(\boldsymbol{\omega}(1)\right)=\underset{\pi}{\max}V_{1}^{\pi}\left(\boldsymbol{\omega}(1)\right), & \textrm{with} & \pi^{\ast}=\arg\max_{\pi}V_{1}^{\pi}\left(\boldsymbol{\omega}(1)\right)\label{eq:optimization_problem}
\end{eqnarray}

\subsection{Transition Probabilities\label{Sec: Transition_probabilities}}

The belief transition probabilities, given decision $\mathcal{\mathcal{U}}(t)=\mathcal{\mathcal{U}}$
and $\boldsymbol{\omega}(t)=\boldsymbol{\omega}=[\omega_{1},...,\omega_{M}]$,
are
\begin{equation}
p_{\boldsymbol{\omega\omega'}}^{(\mathcal{U})}=\Pr\left[\boldsymbol{\omega}(t+1)=\boldsymbol{\omega}'|\boldsymbol{\omega}(t)=\boldsymbol{\omega},\mathcal{\mathcal{U}}(t)=\mathcal{\mathcal{U}}\right]=\prod\limits _{i=1}^{M}\Pr[\omega_{i}(t+1)=\omega_{i}^{\prime}|\omega_{i}(t)=\omega_{i},\mathcal{\mathcal{U}}(t)=\mathcal{\mathcal{U}}],\label{eq:p_ww}
\end{equation}
where $\boldsymbol{\omega}(t+1)=\boldsymbol{\omega}'=[\omega'_{1},...,\omega'_{M}]$,
while the distribution of entry $\omega_{i}(t+1)$ is (see Fig. \ref{fig:Full_DMC})
\begin{equation}
\Pr[\omega_{i}(t+1)=\omega_{i}'|\omega_{i}(t)=\omega_{i},\mathcal{\mathcal{U}}(t)=\mathcal{\mathcal{U}}]=\left\{ \begin{array}{ll}
\omega_{i} & \text{if }\omega_{i}'=p_{11}^{(1)}\textrm{ and }U_{i}\in\mathcal{U}\\
\left(1-\omega_{i}\right) & \text{if }\omega_{i}'=p_{01}^{(1)}\textrm{ and }U_{i}\in\mathcal{U}\\
1 & \text{if }\omega_{i}'=\tau_{0}^{(1)}(\omega_{i})\textrm{ and }U_{i}\notin\mathcal{U}
\end{array}\right.,\label{eqn: belief_update}
\end{equation}
where we have defined the deterministic function
\begin{eqnarray}
\tau_{0}^{(1)}(\omega) & = & \Pr[Q_{i}(t+1)=1|\omega_{i}(t)=\omega,U_{i}\notin\mathcal{U}(t)]=\omega p_{11}^{(0)}+(1-\omega)p_{01}^{(0)}=\omega\delta_{0}+p_{01}^{(0)}\label{eqn: tau1}
\end{eqnarray}
to indicate the next slot's belief when $U_{i}$ is not scheduled
($U_{i}\notin\mathcal{U}(t)$), with $\delta_{0}=p_{11}^{(0)}-p_{01}^{(0)}\geq0$
due to inequalities (\ref{eqn: inequality_transition_probabilities}).
Eq. (\ref{eqn: tau1}) follows from Fig. \ref{fig:Full_DMC}-a), since
the next slot's belief is either $p_{11}^{(0)}$ if $Q_{i}(t)=1$
(w.p. $\omega$) or $p_{01}^{(0)}$ if $Q_{i}(t)=0$ (w.p. $\left(1-\omega\right)$).
A generalization of function $\tau_{0}^{(1)}(\omega)$ that computes
the belief $\omega_{i}(t+k)$ of node $U_{i}$ when it is not scheduled
for $k$ successive slots, e.g., slots $\{t,...,t+k-1\},$ and $\omega_{i}(t)=\omega$,
can be obtained as
\begin{equation}
\tau_{0}^{(k)}(\omega)=\Pr[B_{i}(t+k)=1|\omega_{i}(t)=\omega,U_{i}\notin\mathcal{U}(t),...,U_{i}\notin\mathcal{U}(t+k-1)]=\omega\delta_{0}^{k}+p_{01}^{(0)}\frac{1-\delta_{0}^{k}}{1-\delta_{0}}.\label{eqn: tau_function}
\end{equation}
Eq. (\ref{eqn: tau_function}) can be obtained recursively from (\ref{eqn: tau1})
as $\tau_{0}^{(k)}(\omega)=\tau_{0}^{(1)}(\tau_{0}^{(k-1)}(\omega))$,
for all $k\geq1$, with $\tau_{0}^{(0)}(\omega)=\omega$.

Under assumptions (\ref{eqn: inequality_transition_probabilities}),
it is easy to verify that function (\ref{eqn: tau1}) satisfies the
inequalities
\begin{align}
 & p_{11}^{(1)}\leq p_{01}^{(1)}\leq\tau_{0}^{(1)}(\omega)\text{, for all }0\leq\omega\leq1,\textrm{ and}\label{eqn: properties_tau_functions_line_1}\\
 & \tau_{0}^{(1)}(\omega)\leq\tau_{0}^{(1)}(\omega^{\prime})\text{, for all }\omega\leq\omega^{\prime}\text{ with }0\leq\omega,\omega^{\prime}\leq1.\label{eqn: properties_tau_functions_line_2}
\end{align}
The inequalities in (\ref{eqn: properties_tau_functions_line_1})
guarantee that the belief of a non-scheduled node is always larger
than that of a scheduled one, while the inequality (\ref{eqn: properties_tau_functions_line_2})
says that the belief ordering of two non-scheduled nodes is maintained
across a slot. These inequalities play a crucial role in the analysis
below.

\subsection{Optimality Equations}

The dynamic programming (DP) formulation of problem (\ref{eq:optimization_problem})
(see e.g., \cite{book:Puterman}) allows to express the throughput
recursively over the horizon $\{t,...,T\}$, under policy $\pi$ and
initial belief $\boldsymbol{\omega}$, as 
\begin{align}
V_{t}^{\pi}(\boldsymbol{\omega}) & =\sum_{j=t}^{T}\beta^{j-t}\mathrm{E}^{\mathbf{\pi}}\left[R\left(\boldsymbol{\omega}(j)\mathbf{,}\mathcal{U}^{\mathbf{\pi}}(j)\right)|\boldsymbol{\omega}(t)=\boldsymbol{\omega}\right]=R\left(\boldsymbol{\omega}\mathbf{,}\mathcal{U}^{\mathbf{\pi}}(t)\right)+\beta\sum_{\mathbf{\omega}'}p_{\boldsymbol{\omega\omega}'}^{(\mathcal{U}^{\pi})}V_{t+1}^{\pi}(\boldsymbol{\omega}'),\label{eq: value_DP}
\end{align}
where $V_{t}^{\pi}(\cdot)=0$ for $t>T$. The DP optimality conditions
(or \emph{Bellman equations}) are then expressed in terms of the \textit{value
function} $V_{t}^{\ast}(\boldsymbol{\omega})=\max_{\pi}V_{t}^{\pi}(\boldsymbol{\omega})$,
which represents the optimal throughput in the interval $\{t,...,T\}$,
and it is given by
\begin{eqnarray}
V_{t}^{\ast}(\boldsymbol{\omega}) & = & \underset{\mathcal{U}(t)=\mathcal{U}\subseteq\{U_{1},...,U_{M}\}}{\max}\left\{ R(\mathbf{\boldsymbol{\omega}},\mathcal{U})+\beta\sum_{\mathbf{\omega}'}p_{\boldsymbol{\omega\omega}'}^{(\mathcal{U})}V_{t+1}^{*}(\boldsymbol{\omega}')\right\} .\label{eq:DP-optimality-equation}
\end{eqnarray}
Note that, since the nodes are stochastically equivalent, the value
function (\ref{eq:DP-optimality-equation}) only depends on the numerical
values of the entries of the belief vector $\boldsymbol{\omega}$
regardless of the way it is ordered. Finally, an \emph{optimal policy}
$\pi^{\ast}=\left[\mathcal{U}^{\ast}(1),...,\mathcal{U}^{\ast}(T)\right]$
(see (\ref{eq:optimization_problem})) is such that $\mathcal{U}^{\ast}(t)$
attains the maximum in the condition (\ref{eq:DP-optimality-equation})\emph{
}for $t\in\{1,2,...,T\}$.

\section{Optimality of the Myopic Policy\label{sec:Optimality-of-MP}}

We now define the myopic policy (MP) and show that, under assumptions
(\ref{eq:RR_conditions}), it is a round-robin (RR) policy that schedules
the nodes periodically and that it is optimal for problem (\ref{eq:optimization_problem}).

\subsection{The Myopic Policy is Round-Robin\label{sub:The_MP_Round_robin}}

The MP $\pi^{MP}=\{\mathcal{U}^{MP}(1),...,\mathcal{U}^{MP}(T)\}$,
with throughput $V_{t}^{MP}(\cdot)$, is the greedy policy that schedules
at each slot the $K$ nodes with the largest beliefs so as to maximize
the immediate reward (\ref{eqn: reward_belief_MDP_line_1}), that
is, we have
\begin{align}
\mathcal{U}^{MP}(t) & =\arg\max_{\mathcal{U}}R(\mathbf{\boldsymbol{\omega}}(t),\mathcal{U})=\arg\max_{\mathcal{U}}\sum_{U_{i}\in\mathcal{U}}\omega_{i}(t).\label{eqn: Myopic_policy_maximiz}
\end{align}

\begin{prop}
\emph{\label{Prop: MP_Structure}Under assumptions} \emph{(\ref{eq:RR_conditions}),
the MP (\ref{eqn: Myopic_policy_maximiz}), given an initial belief
$\boldsymbol{\omega}'(1)$, is a RR policy that operates as follows:
}\textbf{\emph{1)}}\emph{ Sort vector $\boldsymbol{\omega}'(1)$ in
a decreasing order to obtain $\boldsymbol{\omega}(1)=\left[\omega_{1}(1),...,\omega_{M}(1)\right]$
such that $\omega_{1}(1)\geq...\geq\omega_{M}(1)$. Re-number the
nodes so that $U_{i}$ has belief $\omega_{i}(1)$; }\textbf{\emph{2)}}\emph{
Divide the nodes into $m$ groups of $K$ nodes each, so that the
$g$th group $\mathcal{G}_{g}$, $g\in\{1,...,m\}$, contains all
nodes $U_{i}$ such that $g=\left\lfloor \frac{i-1}{K}\right\rfloor +1$,
namely: $\mathcal{G}_{1}=\{U_{1},...,U_{K}\},$ $\mathcal{G}_{2}=\{U_{K+1},...,U_{2K}\},$
and so on; }\textbf{\emph{3)}}\emph{ Schedule the groups in a RR fashion
with period $m$ slots, so that groups $\mathcal{G}_{1},...,\mathcal{G}_{m},\mathcal{G}_{1},...$
are sequentially scheduled at slot $t=1,...,m,m+1,...$ and so on. }\end{prop}
\begin{IEEEproof}
According to (\ref{eqn: Myopic_policy_maximiz}), the first scheduled
set is $\mathcal{U}^{MP}(1)$ $=\mathcal{G}_{1}=\{U_{1},U_{2},...,U_{K}\}$.
The beliefs are then updated through (\ref{eqn: belief_update}).
Recalling (\ref{eqn: properties_tau_functions_line_1}), the scheduled
nodes, in $\mathcal{G}_{1}$, have their belief updated to either
$p_{11}^{(1)}$ or $p_{01}^{(1)}$, which are both smaller than the
belief of any non-scheduled node in $\{U_{1},...,U_{M}\}\setminus\mathcal{G}_{1}$.
Moreover, the ordering of the non-scheduled nodes' beliefs is preserved
due to (\ref{eqn: properties_tau_functions_line_2}). Hence, the second
scheduled group is $\mathcal{U}^{MP}(2)$ $=\mathcal{G}_{2}$, the
third is $\mathcal{U}^{MP}(3)$ $=\mathcal{G}_{3},$ and so on. This
proves that the MP, upon an initial ordering of the beliefs, is a
RR policy.
\end{IEEEproof}
We emphasize that the MP sorts the beliefs of the nodes only at the
first slot in which it is operated, and then it keeps scheduling the
groups of nodes according to their initial ordering, without requiring
to recalculate the beliefs.

\subsection{Optimality of the Myopic Policy\label{Sec: optimality_myopic_policy}}

We now prove the optimality of the MP by showing that it satisfies
the Bellman equations\emph{ }(\ref{eq:DP-optimality-equation}). To
start with, let us consider a RR policy $\pi^{RR}$ that operates
according to steps \textbf{2)} and \textbf{3)} of Proposition \ref{Prop: MP_Structure}
(i.e., without re-ordering the initial belief), and let its throughput
(\ref{eq: value_DP}) be denoted by $V_{t}^{RR}(\mathbf{\boldsymbol{\omega}})$.
Note that, when the initial belief $\boldsymbol{\omega}$ is ordered
so that \emph{$\omega_{1}\geq...\geq\omega_{M}$}, then $V_{t}^{RR}(\boldsymbol{\omega})=V_{t}^{MP}(\boldsymbol{\omega})$.
Based on backward induction arguments similarly to \cite{art:Zhao-myopic,art:Ahmad-vector-myopic},
the following lemma establishes a sufficient condition for the optimality
of the MP. 
\begin{lem}
\emph{\label{Lem: optimality_lemma} Assume that the MP is optimal
at slot $t+1,...,T$, i.e., it satisfies (\ref{eq:DP-optimality-equation}).
To show that the MP is optimal also at slot $t$ it is sufficient
to prove the inequality 
\begin{equation}
V_{t}^{RR}(\boldsymbol{\omega}_{\mathcal{S}},\boldsymbol{\omega}_{\mathcal{S}^{c}})\leq V_{t}^{MP}(\boldsymbol{\omega}_{\mathcal{S}},\boldsymbol{\omega}_{\mathcal{S}^{c}})=V_{t}^{RR}(\omega_{1},\omega_{2},...,\omega_{M}),\quad\textrm{for all }\omega_{1}\geq\omega_{2}\geq...\geq\omega_{M}\label{eqn: optimality_lemma_inequality}
\end{equation}
and all sets $\mathcal{S}\subseteq\{1,...,M\}$ of $K$ elements,
with the elements in $\boldsymbol{\omega}_{\mathcal{S}^{c}}$ decreasingly
ordered.}\end{lem}
\begin{IEEEproof}
Since the MP is optimal from $t+1$ onward by assumption, it is sufficient
to show that scheduling $K$ nodes with arbitrary beliefs at slot
$t$ and then following the MP from slot $t+1$ on is no better than
following the MP immediately at slot $t$. The performance of the
former policy is given by the left-hand side (LHS) of (\ref{eqn: optimality_lemma_inequality}).
In fact $V_{t}^{RR}(\boldsymbol{\omega}_{\mathcal{S}},\boldsymbol{\omega}_{\mathcal{S}^{c}})$,
for any set $\mathcal{S}$, represents the throughput of a policy
that schedules the $K$ nodes with beliefs $\boldsymbol{\omega}_{\mathcal{S}}$
at slot $t$, and then operates as the MP from $t+1$ onward, since
beliefs in \emph{$\boldsymbol{\omega}_{\mathcal{S}^{c}}$} are decreasingly
ordered. The MP's performance is instead given by the right-hand side
(RHS) of (\ref{eqn: optimality_lemma_inequality}). Note that, for
$t=T$, it is immediate to verify that the MP is optimal. This concludes
the proof.\end{IEEEproof}
\begin{thm}
\emph{\label{Thr: Theorem_1}Under assumptions (\ref{eq:RR_conditions})
the MP is optimal for problem (\ref{eq:optimization_problem}), so
that $\pi^{MP}=\pi^{\ast}$.}\end{thm}
\begin{IEEEproof}
To start with, we first prove in Appendix \ref{App:Proof-of-Theorem}
that the inequality \emph{
\begin{eqnarray}
V_{t}^{RR}(\omega_{1},...,\omega_{j},y,x,...,\omega_{M}) & \leq & V_{t}^{RR}(\omega_{1},...,\omega_{j},x,y,...,\omega_{M})\label{eqn: inequality_x_and_y}
\end{eqnarray}
}holds for any $x\geq y$, with $0\leq j\leq M-2$, and for all $t\in\{1,...,T\}$
and beliefs $\omega_{k}$ (not necessarily ordered), with $k\in\{1,...,M\}$.
Inequality (\ref{eqn: inequality_x_and_y}) for $j=0$ is intended
as $V_{t}^{RR}(y,x,...,\omega_{M})\leq V_{t}^{RR}(x,y,...,\omega_{M})$.
If (\ref{eqn: inequality_x_and_y}) holds, then inequality (\ref{eqn: optimality_lemma_inequality})
is satisfied for all $\omega_{1}\geq...\geq\omega_{M}$ and all subsets
$\mathcal{S}\subseteq\{1,...,M\}$ of $K$ elements. In fact, (\ref{eqn: inequality_x_and_y})
states that the throughput of the RR policy never increases when,
for any pair of adjacent nodes, the one with the smallest belief of
the pair is scheduled first. Hence, by starting from the RHS of (\ref{eqn: optimality_lemma_inequality})
(i.e., $V_{t}^{RR}(\omega_{1},\omega_{2},...,\omega_{M})$) and by
applying a convenient number of successive switchings between pair
of adjacent elements of vector $[\omega_{1},\omega_{2},...,\omega_{M}]$
to achieve \emph{$[\boldsymbol{\omega}_{\mathcal{S}},\boldsymbol{\omega}_{\mathcal{S}^{c}}]$},
for any \emph{$\mathcal{S}$}, we can obtain a cascade of inequalities
through (\ref{eqn: inequality_x_and_y}) (one for each switching),
which guarantees that (\ref{eqn: optimality_lemma_inequality}) holds.
By Lemma \ref{Lem: optimality_lemma} this is sufficient to prove
that the MP is optimal, since the inequality (\ref{eqn: optimality_lemma_inequality})
holds for any arbitrary $t$.
\end{IEEEproof}

\subsection{Extension to the Infinite-Horizon Case\label{Sec: infinite_horizon}}

We now briefly describe the extension of problem (\ref{eq:optimization_problem})
to the infinite-horizon case, for which the throughput under policy
$\pi$ and its optimal value are given by (see e.g., \cite{art:Zhao-myopic})
\begin{eqnarray}
V^{\pi}\left(\mathbf{\boldsymbol{\omega}}(1)\right)=\sum_{t=1}^{\infty}\beta^{t-1}\mathrm{E}^{\mathbf{\pi}}\left[R\left(\mathbf{\boldsymbol{\omega}}(t)\mathbf{,}\mathcal{U}^{\pi}(t)\right)|\mathbf{\boldsymbol{\omega}}(1)\right], & \textrm{and} & V^{\ast}\left(\mathbf{\boldsymbol{\omega}}(1)\right)=\max_{\pi}V^{\pi}\left(\mathbf{\boldsymbol{\omega}}(1)\right),\label{eq:discounted_infinite_horizon_pair}
\end{eqnarray}
where the optimal policy is $\pi^{\ast}=\arg\max_{\pi}V^{\pi}\left(\mathbf{\boldsymbol{\omega}}(1)\right)$
and $0\leq\beta<1$. From standard DP theory \cite{book:Puterman},
the optimal policy $\pi^{\ast}$ is \textit{\emph{stationary}}, so
that the optimal decision $\mathcal{U}^{*}(t)$ is a function of the
current state $\mathbf{\omega}(t)$ only, independently of slot $t$
\cite{book:Puterman}. By following the same reasoning as in \cite[Theorem 3]{art:Zhao-myopic},
it can be shown that the optimality of the MP for the finite-horizon
setting implies the optimality also for the infinite-horizon scenario.
Moreover, by following \cite[Theorem 4]{art:Zhao-myopic} it can be
shown that the MP is optimal also for the undiscounted average reward
criterion (i.e., $V_{avg}^{\pi}\left(\mathbf{\boldsymbol{\omega}}(1)\right)=lim_{T\rightarrow\infty}\frac{1}{T}\sum_{t=1}^{\infty}\mathrm{E}^{\mathbf{\pi}}\left[R\left(\mathbf{\boldsymbol{\omega}}(t)\mathbf{,}\mathcal{U}^{\pi}(t)\right)|\mathbf{\boldsymbol{\omega}}(1)\right]$).

\section{Optimality of the Whittle Index Policy\label{sec:Whittle_section}}

In this section, we briefly review the Whittle index policy for RMAB
problems \cite{book:Gittins}, and then focus on the infinite-horizon
scenario of Sec. \ref{Sec: infinite_horizon}, when conditions (\ref{eqn: inequality_transition_probabilities})
are specialized to
\begin{equation}
0=p_{11}^{(1)}\leq p_{01}^{(1)}=p_{01}^{(0)}=p_{01}\leq p_{11}^{(0)}=1,\label{eqn: inequality_transition_probabilities_Whittle}
\end{equation}
and where the task queues are of capacity one. We show that under
the assumption (\ref{eqn: inequality_transition_probabilities_Whittle})
(see Sec. \ref{sub:Related-Work-and} for a discussion on these conditions),
the RMAB at hand is indexable and we calculate its Whittle index in
closed-form. We then show that the Whittle index policy is equivalent
to the MP, and thus optimal for the problem (\ref{eq:discounted_infinite_horizon_pair}).

We emphasize that, our results provide a rare example \cite{book:Gittins}
in which, as in \cite{art:Zhao-whittle}, not only indexability is
established, but also the Whittle index is obtained in closed form
and the Whittle policy proved to be optimal. It is finally remarked
that our proof technique is inspired by \cite{art:Zhao-whittle},
but the different system model poses new challenges that require significant
work.

\subsection{Whittle Index}

The Whittle index policy assigns a numerical value $W(\omega_{i})$
to each state $\omega_{i}$ of node $U_{i}$, referred to as \textit{index},
to measure how rewarding it is to schedule node $U_{i}$ in the current
slot. The $K$ nodes with the largest index are then scheduled in
each slot. As detailed below, the Whittle index is calculated independently
for each node, and thus the Whittle index policy is not generally
optimal for RMAB problems. Moreover, even the existence of a well-defined
Whittle index is not guaranteed \cite{book:Gittins}. To study the
indexability and the Whittle index for the RMAB at hand, we can focus
on a restless single-armed bandit (RSAB) model, as defined below \cite{book:Gittins}.
A RSAB is a RMAB with a single arm, in which the only decision that
needs to be taken by the CC is whether activating the (single) arm
or not (i.e., keep it passive).

\subsubsection{RSAB with Subsidy for Passivity\label{Sec: RSAB}}

The Whittle index is based on the concept of \emph{subsidy for passivity},
whereby the CC is given a subsidy $m\in\mathbb{R}$ when the arm is
not scheduled. At each slot $t$, the CC, based on the state $\omega(t)$
of the arm, can decide to activate (or schedule) it, i.e., to set
$u(t)=1$, obtaining an immediate reward $R_{m}(\omega(t),1)=\omega(t)$.
If, instead, the arm is kept passive, i.e., $u(t)=0$, a reward $R_{m}(\omega(t),0)=m$
equal to the subsidy is accrued. The state $\omega(t)$ evolves through
(\ref{eqn: belief_update}), which under (\ref{eqn: inequality_transition_probabilities_Whittle})
and adapted to the simplified notation used here becomes 
\begin{equation}
\omega(t+1)=\left\{ \begin{array}{llc}
0 & \text{w.p. }\omega(t) & \text{if }u(t)=1\\
p_{01} & \text{w.p. }\left(1-\omega(t)\right) & \text{if }u(t)=1\\
\tau_{0}^{(1)}(\omega(t)) & \text{w.p. }1 & \text{if }u(t)=0
\end{array}\right..\label{eqn: belief_update_SAB}
\end{equation}
 The throughput, given policy $\pi=\left\{ u^{\pi}(1),u^{\pi}(2),...\right\} $
and initial belief $\omega(1)$, is 
\begin{equation}
V_{m}^{\pi}\left(\omega(1)\right)=\sum_{t=1}^{\infty}\beta^{t-1}\mathrm{E}^{\pi}\left[R_{m}(\omega(t),u^{\pi}(t))|\omega(1)\right].\label{eqn: discounted_avg_reward_SAB}
\end{equation}
The \textit{\emph{optimal}} throughput is $V_{m}^{\ast}\left(\omega(1)\right)=\max_{\pi}V_{m}^{\pi}\left(\omega(1)\right)$,
while the optimal policy \newline $\pi^{\ast}=\arg\max_{\pi}V_{m}^{\pi}\left(\omega(1)\right)$
is stationary in the sense that the optimal decisions $u_{m}^{\ast}(\omega)\in\{0,1\}$
are functions of the belief $\omega$ only, independently of slot
$t$ \cite{art:Zhao-whittle}. Removing the slot index from the initial
belief, the optimal throughput $V_{m}^{\ast}\left(\omega\right)$
and the optimal decision $u_{m}^{\ast}(\omega)$ satisfy the following
DP optimality equations for the infinite-horizon scenario (see \cite{art:Zhao-whittle})
\begin{eqnarray}
V_{m}^{\ast}(\omega) & = & \underset{u\in\{0,1\}}{\max}\left\{ V_{m}(\omega|u)\right\} ,\label{eqn: optimality_equation_SAB}\\
\textrm{and }u_{m}^{\ast}(\omega) & = & \arg\underset{u\in\{0,1\}}{\max}\left\{ V_{m}(\omega|u)\right\} .\label{eqn: umo}
\end{eqnarray}
In (\ref{eqn: optimality_equation_SAB})-(\ref{eqn: umo}) we defined
$V_{m}(\omega|u),$ $u\in\{0,1\}$, as the throughput (\ref{eqn: discounted_avg_reward_SAB})
of a policy that takes action $u$ at the current slot and then uses
the optimal policy $u_{m}^{\ast}\left(\omega\right)$ onward, we have
\begin{eqnarray}
V_{m}(\omega|0) & = & m+\beta V_{m}^{\ast}(\tau_{0}^{(1)}\left(\omega\right)),\textrm{ \textrm{and}}\label{eqn: passive_action}\\
V_{m}(\omega|1) & = & \omega+\beta\left[\omega V_{m}^{\ast}(0)+(1-\omega)V_{m}^{\ast}(p_{01})\right].\label{eqn: active_action}
\end{eqnarray}

\subsubsection{Indexability and Whittle Index\label{eqn: Definition_indexability}}

We use the notation of \cite{art:Zhao-whittle} to define indexability
and Whittle index for the RSAB at hand. We first define the so called
\textit{passive set} 
\begin{equation}
\mathcal{P}(m)=\left\{ \omega\text{: }0\leq\omega\leq1\text{ and }u_{m}^{\ast}(\omega)=0\right\} ,\label{eqn: Passive_set}
\end{equation}
as the set that contains all the beliefs $\omega$ for which the passive
action is optimal (i.e., all $0\leq\omega\leq1$ such that $V_{m}(\omega|0)\geq V_{m}(\omega|1)$,
see (\ref{eqn: passive_action})-(\ref{eqn: active_action})) under
the given subsidy for passivity $m\in\mathbb{R}$. The RMAB at hand
is said to be \textit{indexable} if the passive set $\mathcal{P}(m)$,
for the associated RSAB problem%
\footnote{Note that in a RMAB with arms characterized by different statistics
this condition must be checked for all arms.%
}, is monotonically increasing as $m$ increases within the interval
$(-\infty,+\infty),$ in the sense that $\mathcal{P}(m^{\prime})\subseteq$
$\mathcal{P}(m)$ if $m^{\prime}\leq m$ and $\mathcal{P}(-\infty)=\emptyset$
and $\mathcal{P}(+\infty)=[0,1]$.

If the RMAB is indexable, the Whittle index $W(\omega)$ for each
arm with state $\omega$ is the infimum subsidy $m$ such that it
is optimal to make the arm passive. Equivalently, the Whittle index
$W(\omega)$ is the infimum subsidy $m$ that makes passive and active
actions equally rewarding, i.e., 
\begin{eqnarray}
W(\omega) & = & \inf\left\{ m\text{: }u_{m}^{\ast}(\omega)=0\right\} =\inf\left\{ m\text{: }V_{m}\left(\omega|0\right)=V_{m}\left(\omega|1\right)\right\} .\label{eqn: definition_whittle_index_DEF}
\end{eqnarray}

\subsection{Optimality of the Threshold Policy}

\begin{figure}[h!]
\centering \includegraphics[clip,width=15cm]{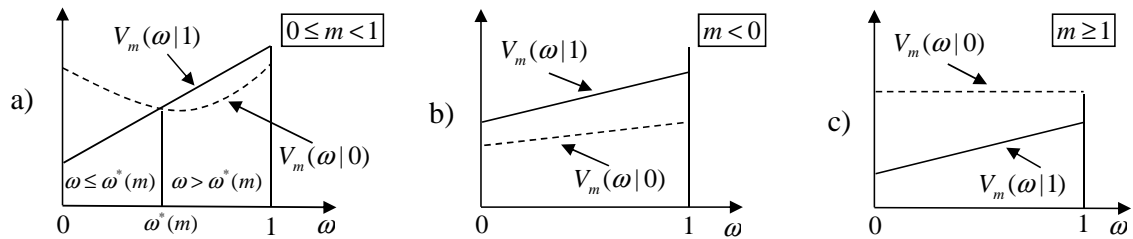}
\caption{Illustration of the optimality of a threshold policy for different
values of the subsidy for passivity $m$: a) $0\leq m<1$; b) $m<0$;
c) $m\geq1$.}

\label{fig:Extrema_inequalities}
\end{figure}

Here, we show that the optimal policy $u_{m}^{\ast}(\omega)$ for
the RSAB of Sec. \ref{Sec: RSAB} is a threshold policy over the belief
$\omega$. This is crucial in our proof of indexability of the RMAB
at hand given in Sec. \ref{sub:Indexability-and-Whittle}. To this
end, we observe that: \textit{i}) function $V_{m}(\omega|1)$ in (\ref{eqn: active_action})
is linear over the belief $\omega$; \textit{ii}) function $V_{m}(\omega|0)=m+\beta V_{m}^{\ast}(\tau_{0}^{(1)}\left(\omega\right))$
in (\ref{eqn: passive_action}) is convex over $\omega$, since the
value function $V_{m}^{\ast}(\omega)$ is convex for the problem at
hand (see \cite{art:Zhao-whittle,art:Cassandra}). We need the following
lemma.
\begin{lem}
\begin{flushleft}
\label{Lem: Extrema_inequalities_lemma}\emph{The following inequalities
hold:}\begin{subequations}
\begin{align}
\mathcal{\textrm{a) For }} & 0\leq m<1: & \mathcal{\textrm{a.1) }}V_{m}(0|1)\leq V_{m}(0|0)\leq V_{m}(1|1);\quad\mathcal{\textrm{a.2) }}V_{m}(1|0)\leq V_{m}(1|1);\label{eqn: extrema_inequality_1_m_0_1}\\
\mathcal{\textrm{b) For }} & m<0: & \mathcal{\textrm{b.1) }}V_{m}(0|0)\leq V_{m}(0|1)\leq V_{m}(1|1);\quad\mathcal{\textrm{b.2) }}V_{m}(1|0)\leq V_{m}(1|1);\label{eqn: extrema_inequality_1_m_negative}\\
\mathcal{\textrm{c) For }} & m\geq1: & \mathcal{\textrm{c.1) }}V_{m}(0|0)\leq V_{m}(1|1)\leq V_{m}(0|1);\quad\mathcal{\textrm{c.2) }}V_{m}(1|1)\leq V_{m}(1|0).\label{eqn: extrema_inequality_1_m_bigger_1}
\end{align}
\end{subequations}
\par\end{flushleft}\end{lem}
\begin{IEEEproof}
See Appendix \ref{App: proof_extrema_inequality_m_0_1}.
\end{IEEEproof}
Leveraging Lemma \ref{Lem: Extrema_inequalities_lemma}, we can now
establish the optimality of a threshold policy $u_{m}^{\ast}(\omega)$.
\begin{prop}
\label{Prop: Threshold_policy}\emph{The optimal policy $u_{m}^{\ast}(\omega)$
in (\ref{eqn: umo}) for subsidy $m\in\mathbb{R}$ is given by}

\begin{equation}
u_{m}^{\ast}(\omega)=\left\{ \begin{array}{cc}
1, & \text{if }\omega>\omega^{\ast}(m)\\
0, & \text{if }\omega\leq\omega^{\ast}(m)
\end{array}\right.,\label{eqn: optimal_threshold_policy}
\end{equation}
\emph{where $\omega^{\ast}(m)\in\mathbb{R}$ is the optimal threshold
for a given subsidy $m$. The optimal threshold $\omega^{\ast}(m)$
is $0\leq\omega^{\ast}(m)\leq1$ if $0\leq m<1$, while it is arbitrary
negative for $m<0$ and arbitrary greater than unity for $m\geq1$.
In other words we have $u_{m}^{\ast}(\omega)=1$ if $m<0$ and $u_{m}^{\ast}(\omega)=0$
if $m\geq1$. }\end{prop}
\begin{IEEEproof}
We start by showing that (\ref{eqn: optimal_threshold_policy}), for
$0\leq m<1$, satisfies (\ref{eqn: umo}) and is thus an optimal policy.
To see this, we refer to Fig. \ref{fig:Extrema_inequalities}, where
we sketch functions $V_{m}(\omega|1)$ and $V_{m}(\omega|0)$ for
different values of the subsidy $m$. From (\ref{eqn: umo}), we have
that $u_{m}^{\ast}(\omega)=1$ for all $\omega$ such that $V_{m}(\omega|1)>V_{m}(\omega|0)$
and $u_{m}^{\ast}(\omega)=0$ otherwise. For $0\leq m<1$, from the
inequalities of Lemma \ref{Lem: Extrema_inequalities_lemma}-a), the
linearity of $V_{m}(\omega|1)$ and the convexity of $V_{m}(\omega|0)$,
it follows that there is only one intersection $\omega^{\ast}(m)$
between $V_{m}(\omega|1)$ and $V_{m}(\omega|0)$ with $0\leq\omega^{\ast}(m)\leq1$,
as shown in Fig. \ref{fig:Extrema_inequalities}-a). Instead, when
$m<0$, by Lemma \ref{Lem: Extrema_inequalities_lemma}-b), arm activation
is always optimal, that is, $u_{m}^{\ast}(\omega)=1$, since $V_{m}(\omega|1)>V_{m}(\omega|0)$
for any $0\leq\omega\leq1$ as shown in Fig. \ref{fig:Extrema_inequalities}-b).
Conversely, when $m\geq1$, by Lemma \ref{Lem: Extrema_inequalities_lemma}-c),
it follows that passivity is always optimal, that is, $u_{m}^{\ast}(\omega)=0$,
since $V_{m}(\omega|0)\geq V_{m}(\omega|1)$ for any $0\leq\omega\leq1$
as shown in Fig. \ref{fig:Extrema_inequalities}-c).
\end{IEEEproof}

\subsection{Closed-Form Expression of the Value Function}

By leveraging the optimality of the threshold policy (\ref{eqn: optimal_threshold_policy})
we derive a closed-form expression of $V_{m}^{\ast}(\omega)$ in (\ref{eqn: optimality_equation_SAB}),
being a key step in establishing the RMAB's indexability in Sec. \ref{sub:Indexability-and-Whittle}.

Notice that function $\tau_{0}^{(k)}(\omega)$ in (\ref{eqn: tau_function}),
when specialized to conditions (\ref{eqn: inequality_transition_probabilities_Whittle}),
becomes 
\begin{equation}
\tau_{0}^{(k)}(\omega)=1-(1-p_{01})^{k}(1-\omega),\label{eqn: tau_h_whittle}
\end{equation}
 which is a monotonically increasing function of $k$, so that $\tau_{0}^{(k)}(\omega)\geq\tau_{0}^{(i)}(\omega)$
for any $k\geq i$. Based on such monotonicity, we can define the
average number $L(\omega,\omega^{\prime})$ of slots it takes for
the belief to become larger than $\omega^{\prime}$ when starting
from $\omega$ while the arm is kept passive, as 
\begin{eqnarray}
L(\omega,\omega^{\prime}) & =\min\left\{ k\text{: }\tau_{0}^{\left(k\right)}(\omega)>\omega^{\prime}\right\}  & =\left\{ \begin{array}{lc}
0 & \omega>\omega^{\prime}\\
\left\lfloor \frac{\ln\left(\frac{1-\omega^{\prime}}{1-\omega}\right)}{\ln\left(1-p_{01}\right)}\right\rfloor +1 & \omega\leq\omega^{\prime}\\
\infty & \omega\leq1\leq\omega^{\prime}
\end{array}\right..\label{eqn: definition_expected_crossing_time}
\end{eqnarray}
From (\ref{eqn: definition_expected_crossing_time}) we have $L(\omega,\omega^{\prime})=1$
for $\omega=\omega^{\prime}$ since, without loss of optimality, we
assumed that the passive action is optimal (i.e., $u_{m}^{\ast}(\omega)=0$)
when $V_{m}(\omega|0)=V_{m}(\omega|1)$. For $\omega^{\prime}\geq1$
instead (according to Proposition \ref{Prop: Threshold_policy}),
the arm is always kept passive and thus $L(\omega,\omega^{\prime})=\infty$.
\begin{lem}
\label{Lem: Value_function_closed}\emph{The optimal throughput $V_{m}^{\ast}(\omega)$
in (\ref{eqn: optimality_equation_SAB}) can be written as}
\begin{equation}
V_{m}^{\ast}(\omega)=\frac{1-\beta^{L(\omega,\omega^{\ast}(m))}}{1-\beta}m+\beta^{L(\omega,\omega^{\ast}(m))}V_{m}(\tau_{0}^{\left(L(\omega,\omega^{\ast}(m))\right)}(\omega)|1),\label{eqn: Value_function_closed_form}
\end{equation}
\emph{where $\omega^{\ast}(m)$ is the optimal threshold obtained
from Proposition \ref{Prop: Threshold_policy}.}\end{lem}
\begin{IEEEproof}
According to Proposition \ref{Prop: Threshold_policy}, the optimal
policy $u_{m}^{\ast}(\omega)$ keeps the arm passive as long as the
current belief is $\omega\leq\omega^{\ast}(m)$. Therefore, the arm
is kept passive for $L(\omega,\omega^{\ast}(m))$ slots, during which
a reward $R_{m}(\omega,0)=m$ is accrued in each slot. This leads
to a total reward within the passivity time given by the following
geometric series $\sum_{k=0}^{L(\omega,\omega^{\ast}(m))-1}\beta^{k}m=\frac{1-\beta^{L(\omega,\omega^{\ast}(m))}}{1-\beta}m$,
which corresponds to the first term in the RHS of (\ref{eqn: Value_function_closed_form}).
After $L(\omega,\omega^{\ast}(m))$ slots of passivity, the belief
becomes larger than the threshold $\omega^{\ast}(m)$ and the arm
is activated. The contribution to the value function $V(\omega)$
thus becomes $\beta^{L(\omega,\omega^{\ast}(m))}V_{m}(\tau_{0}^{\left(L(\omega,\omega^{\ast}(m))\right)}(\omega)|1)$,
which is the second term in the RHS of (\ref{eqn: Value_function_closed_form}).
Note that, when $\omega>\omega^{\ast}(m)$, activation is optimal,
and $V^{\ast}(\omega)=V(\omega|1)$.
\end{IEEEproof}
To evaluate $V_{m}^{\ast}(\omega)$ from (\ref{eqn: Value_function_closed_form}),
we only need to calculate $V_{m}(\omega|1)$ since the other terms,
thanks to (\ref{eqn: definition_expected_crossing_time}) are explicitly
given once $\omega^{\ast}(m)$ is obtained from Proposition \ref{Prop: Threshold_policy}.
However, from (\ref{eqn: active_action}), evaluating $V_{m}(\omega|1)$
only requires $V_{m}^{\ast}(0)$ and $V_{m}^{\ast}(p_{01})$, which
are calculated in the lemma below.
\begin{lem}
\label{Lem: optimal_V_0_and_V_h}\emph{We have}\begin{subequations}
\label{eqn: V_0_and_V_h_closed_form}
\begin{align}
 & V_{m}^{\ast}(0)=\frac{\left(m-2m\beta^{L_{m}^{\ast}}+\beta^{L_{m}^{\ast}}\upsilon_{m}^{\ast}-\beta^{L_{m}^{\ast}+1}\upsilon_{m}^{\ast}+m\beta^{L_{m}^{\ast}+1}+m\beta^{L_{m}^{\ast}}\upsilon_{m}^{\ast}-m\beta^{L_{m}^{\ast}+1}\upsilon_{m}^{\ast}\right)}{\left(\beta-1\right)\left(\beta^{L_{m}^{\ast}}-\beta^{L_{m}^{\ast}}\upsilon_{m}^{\ast}+\beta^{L_{m}^{\ast}+1}\upsilon_{m}^{\ast}-1\right)}\label{eqn: V_0_closed_form}\\
 & V_{m}^{\ast}(p_{01})=\frac{\left(m\beta-m\beta^{L_{m}^{\ast}}+\beta^{L_{m}^{\ast}}\upsilon_{m}^{\ast}-\beta^{L_{m}^{\ast}+1}\upsilon_{m}^{\ast}+m\beta^{L_{m}^{\ast}+1}\upsilon_{m}^{\ast}-m\beta^{L_{m}^{\ast}+2}\upsilon_{m}^{\ast}\right)}{\beta\left(\beta-1\right)\left(\beta^{L_{m}^{\ast}}-\beta^{L_{m}^{\ast}}\upsilon_{m}^{\ast}+\beta^{L_{m}^{\ast}+1}\upsilon_{m}^{\ast}-1\right)}\label{eqn: V_h_closed_form}
\end{align}
 \end{subequations} \emph{where we have defined $L_{m}^{\ast}=L(0,\omega^{\ast}(m))$
and $\upsilon_{m}^{\ast}=\tau_{0}^{\left(L(0,\omega^{\ast}(m))\right)}(0).$}\end{lem}
\begin{IEEEproof}
By plugging (\ref{eqn: active_action}) into (\ref{eqn: Value_function_closed_form}),
and evaluating (\ref{eqn: Value_function_closed_form}) for $\omega=0$
and $\omega=p_{01},$ we get a linear system in the two unknowns $V_{m}^{\ast}(0)$
and $V_{m}^{\ast}(p_{01})$, which can be solved leading to (\ref{eqn: V_0_and_V_h_closed_form}).
\end{IEEEproof}

\subsection{Indexability and Whittle Index\label{sub:Indexability-and-Whittle}}

Here, we prove that the RMAB at hand is indexable, we derive the Whittle
index in closed form and show that it is equivalent to the MP and
thus optimal for the RMAB problem (\ref{eq:discounted_infinite_horizon_pair}).
\begin{thm}
\emph{\label{thm:indexability}}\textbf{\emph{a}}\emph{) The RMAB
at hand is indexable and }\textbf{\emph{b}}\emph{) its Whittle index
is}
\begin{equation}
W(\omega)=\frac{\left(1-\beta^{L(0,\omega)}\left(1-\beta\tau_{0}^{L(0,\omega)}(0)\left(1-\beta\right)\left(1-h\right)\right)\right)\omega+\beta^{L(0,\omega)}\tau_{0}^{L(0,\omega)}(0)\left(1-\beta\right)\left(h\beta+1\right)}{(\beta-1)\left(\beta^{L(0,\omega)}\left(1-\beta(1-h)\right)\omega-\left(1+\beta^{L(0,\omega)}\left(\tau_{0}^{L(0,\omega)}(0)(1-\beta)+h\beta\right)\right)\right)}.\label{eqn: Whittle_index}
\end{equation}
\end{thm}
\begin{IEEEproof}
Part \textbf{a}). See Appendix \ref{Sec: Proof_indexability}. Part
\textbf{b}). By (\ref{eqn: definition_whittle_index_DEF}), the Whittle
index $W(\omega)$ of state $\omega$ is the value of the subsidy
$m$ for which activating or not the arm is equally rewarding so that
$V_{m}\left(\omega|0\right)=V_{m}\left(\omega|1\right)$. By using
(\ref{eqn: passive_action})-(\ref{eqn: active_action}) this becomes
$\omega+\beta\left[\omega V_{m}^{\ast}(0)+(1-\omega)V_{m}^{\ast}(p_{01})\right]=m+\beta V_{m}^{\ast}(\tau_{0}^{(1)}\left(\omega\right)).$
Moreover, since the threshold policy is optimal and $\tau_{0}^{(1)}\left(\omega\right)>\omega$,
it follows that, when the belief becomes $\tau_{0}^{(1)}\left(\omega\right)$,
it is optimal to activate the arm and thus $V_{m}^{\ast}(\tau_{0}^{(1)}\left(\omega\right))$
$=V_{m}(\tau_{0}^{(1)}\left(\omega\right)|1)=\beta\tau_{0}^{(1)}(\omega)V_{m}^{\ast}(0)+\beta(1-\tau_{0}^{(1)}(\omega))V_{m}^{\ast}(p_{01})$.
Plugging this result into $V_{m}\left(\omega|0\right)=V_{m}\left(\omega|1\right)$,
along with (\ref{eqn: V_0_closed_form}) and (\ref{eqn: V_h_closed_form}),
leads to (\ref{eqn: Whittle_index}), which concludes the proof.
\end{IEEEproof}
It can be show that the Whittle index $W(\omega)$ in (\ref{eqn: Whittle_index})
is an increasing function of $\omega.$ Therefore, since the Whittle
policy selects the $K$ arms with the largest index at each slot,
we have:
\begin{cor}
\emph{The Whittle index policy is equivalent to the MP and is thus
optimal.}
\end{cor}

\section{Extension to Task Queues of Arbitrary Capacity $C>1$\label{sec:Extension_capacity}}

The problem of characterizing the optimal policies when $C>1$ is
significantly more complicated than for $C=1$ and is left open by
this work. Moreover, since the dimension of the state space of the
belief MDP grows with $C$, even the numerical computation of the
optimal policies is quite cumbersome. Due to these difficulties, here
we compare the performance of the MP, inspired by its optimality for
$C=1$, with a performance upper bound obtained following the relaxation
approach of \cite{art:whittle}. 

\begin{figure}[h!]
\centering \includegraphics[clip,width=10cm]{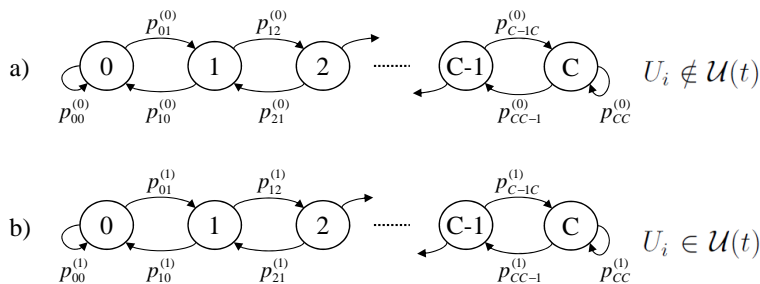}
\caption{Markov model for the evolution of the queue $Q_{i}(t)$, of arbitrary
capacity $C$, when the node $U_{i}$: a) is not scheduled in slot
$t$ (i.e., $U_{i}\notin\mathcal{U}(t)$); b) is scheduled in slot
$t$ (i.e., $U_{i}\in\mathcal{U}(t)$).}

\label{fig:extended_battery}
\end{figure}

\subsection{System Model and Myopic Policy\label{sub:System-Model_extended_capacity}}

Each node $U_{i}$ has a task queue $Q_{i}(t)\in\{0,1,...,C\}$ of
capacity $C$. We consider the Markov model of Fig. \ref{fig:extended_battery}
for the task generation and expiration processes at each node (cf.
Sec. \ref{sub:Markov-Formulation}). The transition probabilities
between queue states when node $U_{i}$ is not scheduled are $p_{xy}^{(0)}=\textrm{Pr}[Q_{i}(t+1)=y|Q_{i}(t)=x,U_{i}\notin\mathcal{U}(t)]$,
whereas when $U_{i}$ is scheduled we have $p_{xy}^{(1)}=\textrm{Pr}[Q_{i}(t+1)=y|Q_{i}(t)=x,U_{i}\in\mathcal{U}(t)]$,
for $x,y\in\{0,1,...,C\}$. When node\emph{ $U_{i}$} is scheduled
at slot $t$, and $Q_{i}(t)\geq1$, one of its task is executed and
it also informs the CC about the number of tasks left in the queue
(observation). We assume that at most one task can be generated (or
dropped) in a slot, so that $p_{xy}^{(u)}=0$ for $y<x-1$ and $y>x+1$,
with $u\in\{0,1\}$ as shown in Fig. \ref{fig:extended_battery}. 

The belief of each $i$th node is represented by a $\left(C\times1\right)$
vector $\boldsymbol{\omega}_{i}=\left[\omega_{i,0},...,\omega_{i,C-1}\right]$
whose $k$th entry $\omega_{i,k}$, for $k\in\{0,1,...,C-1\}$, is
given by (cf. (\ref{eqn: belief_vector_entries})) $\omega_{i,k}=\textrm{Pr}\left[Q_{i}(t)=k|\mathcal{H}(t)\right].$
The immediate reward (\ref{eqn: reward_belief_MDP_line_1}), given
the initial belief vectors $\boldsymbol{\omega}_{1}(t),...,\boldsymbol{\omega}_{M}(t)$
and action $\mathcal{U}$, becomes
\begin{eqnarray}
R(\boldsymbol{\omega}_{1}(t),...,\boldsymbol{\omega}_{M}(t),\mathcal{U}) & = & \sum_{i=1}^{M}\textrm{Pr}\left[Q_{i}(t)>0|\mathcal{H}(t)\right]1(U_{i}\in\mathcal{U})=K-\sum_{i\in\mathcal{U}}\omega_{i,0}(t).\label{eq:AVG_immediate_reward_extended_capacity}
\end{eqnarray}
The performance of interest is the infinite-horizon throughput (\ref{eq:discounted_infinite_horizon_pair}).

\subsubsection{Myopic Policy\label{sub:Myopic_Policy_extended}}

The MP (\ref{eqn: Myopic_policy_maximiz}), specialized to the immediate
reward (\ref{eq:AVG_immediate_reward_extended_capacity}), becomes
\begin{eqnarray}
\mathcal{U}^{MP}(t) & = & \textrm{arg}\underset{\mathcal{U}}{\textrm{max}}R(\boldsymbol{\omega}_{1}(t),...,\boldsymbol{\omega}_{M}(t),\mathcal{U})=\textrm{arg}\underset{\mathcal{U}}{\textrm{min}}\sum_{i\in\mathcal{U}}\omega_{i,0}(t).\label{eq:Myopic_policy_exteded_capacity}
\end{eqnarray}
Note that, unlike Sec. \ref{sub:The_MP_Round_robin}, when $C>1$
the MP does not generally have a RR structure.

\subsection{Upper Bound \label{sub:Whittle_UB}}

Here we derive an upper bound to the throughput (\ref{eq:discounted_infinite_horizon_pair})
by following the approach for general RMAB problems proposed in \cite{art:whittle}.
The upper bound relaxes the constraint that exactly $K$ nodes must
be scheduled in each slot. Specifically, it allows a variable number
$K^{\pi}(t)$ of scheduled nodes in each $t$th slot under policy
$\pi$, with the only constraint that its discounted average satisfies
\begin{equation}
E^{\pi}\left[\sum_{t=1}^{\infty}\beta^{t-1}K^{\pi}(t)\right]=\frac{K}{1-\beta}.\label{eq:AVG_active_users_constraint}
\end{equation}
The advantage of this relaxed version of the scheduling problem is
that it can be tackled by focusing on each single arm independently
from the others \cite{art:whittle,art:Mora}. This is because, by
the symmetry of the nodes, the constraint (\ref{eq:AVG_active_users_constraint})
can be equivalently handled by imposing that each node is active on
average for a discounted time $E^{\pi}[\sum_{t=1}^{\infty}\beta^{t-1}1(U_{i}\in\mathcal{U}^{\pi}(t))]=\frac{K}{M(1-\beta)}$.
We can thus calculate the optimal solution of the relaxed problem
by solving a single RSAB problem. 

We now elaborate on such a RSAB by dropping the node index. Here,
the immediate reward when the arm is in state $\boldsymbol{\omega}$
(a vector since $C>1$, see Sec. \ref{sub:System-Model_extended_capacity}),
and action $u\in\{0,1\}$ is chosen, is $R(\boldsymbol{\omega},u)=1-\omega_{0}$
if $u=1$ and $R(\boldsymbol{\omega},u)=0$ if $u=0$, while the Markov
evolution of the belief follows from Fig. \ref{fig:extended_battery}
and similarly to Sec. \ref{sub:Markov-Formulation}. The problem consists
in optimizing the throughput under the constraint $E^{\pi}[\sum_{t=1}^{\infty}\beta^{t-1}1(U_{i}\in\mathcal{U}^{\pi}(t))]=\sum_{t=1}^{\infty}\beta^{t-1}E^{\pi}[u^{\pi}(t)]=K/(M(1-\beta))$,
as introduced above. Under the assumption that the state $\boldsymbol{\omega}$
belongs to a finite state space $\mathcal{W}$ (to be discussed below),
this optimization can be done by resorting to a linear programming
(LP) formulation \cite{art:Mora}. Specifically, let $z_{\boldsymbol{\omega}}^{(u)}$
be the probability of being in state $\boldsymbol{\omega}$ and selecting
action $u\in\{0,1\}$ under a given policy. The optimization at hand
leads to the following LP \begin{subequations} 
\begin{eqnarray}
\textrm{maximize} & \sum_{\boldsymbol{\omega},u} & R(\boldsymbol{\omega},u)z_{\boldsymbol{\omega}}^{(u)},\label{eq:main_LP}\\
\textrm{subject to}:\;\sum_{\boldsymbol{\omega},u}z_{\boldsymbol{\omega}}^{(u)} & = & 1,\label{eq:LP_probability_constraint}\\
\sum_{\boldsymbol{\omega}}z_{\boldsymbol{\omega}}^{(1)} & = & \frac{K}{M(1-\beta)},\label{eq:LP_AVG_activation}\\
z_{\boldsymbol{\omega}}^{(0)}+z_{\boldsymbol{\omega}}^{(1)} & = & \delta\left(\boldsymbol{\omega}-\boldsymbol{\omega}(1)\right)+\beta\sum_{\boldsymbol{\omega}',u}z_{\boldsymbol{\omega'}}^{(u)}p_{\boldsymbol{\omega}\boldsymbol{\omega}'}^{(u)}\textrm{, for all }\boldsymbol{\omega\in}\mathcal{W},\label{eq:LP_steady_state}
\end{eqnarray}
\label{eq:LP}\end{subequations}where (\ref{eq:LP_AVG_activation})
is the constraint on the average time in which the node is scheduled,
while (\ref{eq:LP_steady_state}) guarantees that $z_{\boldsymbol{\omega}}^{(u)}$
is the stationary distribution \cite{art:Mora}, in which $\delta\left(\boldsymbol{\omega}-\boldsymbol{\omega}(1)\right)=1$
if $\boldsymbol{\omega}=\boldsymbol{\omega}(1)$ and $\delta\left(\boldsymbol{\omega}-\boldsymbol{\omega}(1)\right)=0$
if $\boldsymbol{\omega}\neq\boldsymbol{\omega}(1)$ . Note that, as
discussed in Sec. \ref{sec:System-Model}, the term $p_{\boldsymbol{\omega}\boldsymbol{\omega}'}^{(u)}$
is the probability that the next state is $\boldsymbol{\omega}'$
given that action $u$ is taken in state $\boldsymbol{\omega}$.

We are left to discuss the cardinality of the set $\mathcal{W}$.
While the belief $\boldsymbol{\omega}$ can generally assume any value
in the $C$-dimensional probability simplex, the number of states
actually assumed by $\omega$ during any \emph{limited} horizon of
time is finite due to the finiteness of the action space \cite{art:Cassandra}.
In our problem, since the time horizon is unlimited, this fact alone
is not sufficient to conclude that the set $\mathcal{W}$ is finite.
However, after each $t$th slot in which the arm is activated, the
belief at the $(t+1)$th slot can only takes $C$ values given that
the queue state is learned by the CC. Therefore, the evolution of
the belief is reset after each activation, and in practice, the time
between two activations is finite since the node must be kept active
for a discounted fraction of time $K/\left(M(1-\beta\right)$. Hence,
by constraining the maximum time interval between two activations
to a sufficiently large value, the state space $\mathcal{W}$ remains
finite and the optimal performance is not affected. We used this approach
for the numerical evaluation of the upper bound in Sec. \ref{sub:Numerical-Results}.

\subsection{Numerical Results\label{sub:Numerical-Results}}

We now present some numerical results to compare the performance of
the MP with the upper bound of Sec. \ref{sub:Whittle_UB}. The performance
is the throughput (\ref{eq:discounted_infinite_horizon_pair}) normalized
by its ideal value $K/\left(1-\beta\right)$ that is obtained if the
nodes always have a task to be completed when scheduled.

\begin{figure}[h!]
\centering \includegraphics[clip,width=10cm]{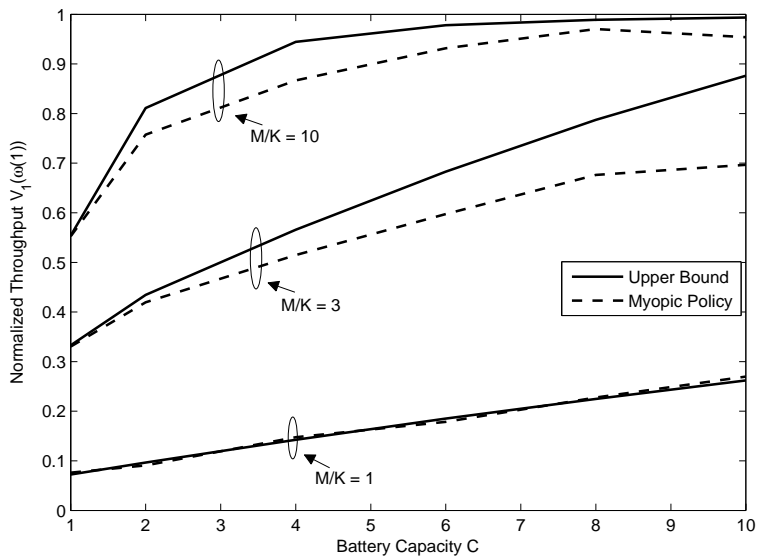}
\caption{Normalized optimal throughput of the MP in (\ref{eq:Myopic_policy_exteded_capacity})
as compared to the upper bound versus the queue capacity $C$ for
different ratios $M/K\in\{1,3,10\}$ (system parameters are $K=3$,
$\beta=0.95$, $\omega_{i,k}(1)=1/(C+1)$ for all $i,k$, $p_{01}^{(0)}=0.15$,
$p_{01}^{(1)}=0.05$, $p_{CC}^{(0)}=0.9$, $p_{CC}^{(1)}=0.05$, $p_{kk-1}^{(0)}=0.05$,
$p_{kk-1}^{(1)}=0.95$, $p_{kk+1}^{(0)}=0.1$, $p_{kk+1}^{(1)}=0$,
for $k\in\{1,C-1\}$).}

\label{fig:AVG_reward_vs_capacity}
\end{figure}

In Fig. \ref{fig:AVG_reward_vs_capacity} we show the normalized throughput
versus the queue capacity $C$ for different ratio $M/K$ between
the number $M$ of nodes and the number $K$ of nodes scheduled in
each slot. We keep $K=3$ fixed and vary $M$. We assume a uniform
distribution for the initial number of tasks in the queues $Q_{i}(1)$
for all the nodes, so that $\omega_{i,k}(1)=1/\left(C+1\right)$ for
all $i,k$. The probabilities that a new task is generated when the
arm is kept passive are $p_{01}^{(0)}=0.15$ and $p_{kk+1}^{(0)}=0.1$,
for $k\in\{1,C-1\}$, while under activation they are $p_{01}^{(1)}=0.05$
and $p_{kk+1}^{(1)}=0$. The probability that a task expires when
the arm is kept passive and activated are $p_{kk-1}^{(0)}=0.05$ and
$p_{kk-1}^{(1)}=0.95$ respectively. The remaining transitions probabilities
are $p_{CC}^{(0)}=0.9$, $p_{CC}^{(1)}=0.05$, while $\beta=0.95$.

From Fig. \ref{fig:AVG_reward_vs_capacity} it can be seen that when
$C$ and/or $M/K$ are small the MP's performance is close to the
upper bound. In fact, for small $M/K$, most of the nodes are scheduled
in each slot and the relaxed system in Sec. \ref{sub:Whittle_UB}
approaches the original one, while for small $C$ we get closer to
the optimality of the MP for $C=1$. For moderate to large values
of $M/K$ and/or $C$ instead, the more flexibility in the relaxed
system enables larger gains over the MP. 

\section{Conclusions}

This paper considers a centralized scheduling problem for independent,
symmetric and time-sensitive tasks under resources constraints. The
problem is to assign a finite number of resources to a larger number
of \emph{nodes }that may have tasks to be completed in their \emph{task
queue}. It is assumed that the \emph{central controller }has no direct
access to the queue of each node. Based on a Markovian modeling of
the task generation and expiration processes, the scheduling problem
is formulated as a partially observable Markov decision process (POMDP)
and then cast into the framework of restless multi-armed bandit (RMAB)
problems. Under the assumption that the task queues are of capacity
one, a greedy, or \emph{myopic policy} (MP), operating in the space
of the a posteriori probabilities (beliefs) of the number of tasks
in the queues, is proved to be optimal, under appropriate assumptions,
for both finite and infinite-horizon throughput criteria. The MP selects
at each slot the nodes with the largest probability of having a task
to be completed. It is shown that the MP is round-robin since it schedules
the nodes periodically. We have also established that the RMAB problem
at hand is indexable, derived the Whittle index in closed form and
shown that the Whittle index policy is equivalent to the MP and thus
it is optimal.

Systems in which the task queues have arbitrary capacities have been
investigated as well by comparing the performance of the MP, which
is generally suboptimal, with an upper bound based on a relaxation
of the scheduling constraint.

Overall, this paper proposes a general framework for resource allocation
that finds applications in several areas of current interest including
communication networks and distributed computing.

\appendices{}

\section{Proof of Theorem \ref{Thr: Theorem_1}\label{App:Proof-of-Theorem}}

The proof is divided into two steps. In the first step we derive the
throughput of the RR policy in closed form, and then we show that
inequality (\ref{eqn: inequality_x_and_y}) holds. 

As for the first step, the throughput for the RR policy (and thus
of the MP) can be calculated as the sum of the contribution of each
node separately (due to the round robin structure). To elaborate,
let us focus on node $U_{i}$, with initial belief $\omega_{i}(1)$,
and assume that $U_{i}\in\mathcal{G}_{1}$. Nodes in group $\mathcal{G}_{1}$
are scheduled at slots $t\in\{1+(j-1)m\}$, for $j\in\{1,2,...\}$.
Let $r_{j}(\omega_{i}(1))=\mathrm{E^{RR}}[\omega_{i}(1+(j-1)m)|\omega_{i}(1)]$
be the average reward accrued by the CC from node $U_{i}$ only, when
scheduling it for the $j$th time at slot $t=1+(j-1)m$ (see the RHS
of (\ref{eqn: reward_belief_MDP_line_1})) (i.e., when operating the
RR policy). At slot $t=1$ we have $r_{1}(\omega_{i}(1))=\omega_{i}(1)$.
To calculate $r_{2}(\omega_{i}(1))$ we first derive the average value
of the belief (see (\ref{eqn: belief_update})) after the slot of
activity in $t=1$ as $\mathrm{E^{RR}}[\omega_{i}(2)|\omega_{i}(1)]=\tau_{1}^{(1)}(\omega_{i}(1))$,
where $\tau_{1}^{(1)}=\omega\delta_{1}+p_{01}^{(1)}$ with $\delta_{u}=\left(p_{11}^{(u)}-p_{01}^{(u)}\right)$
(cf. (\ref{eqn: tau1})). We then account for the $(m-1)$ slots of
passivity by exploiting (\ref{eqn: tau1}), so that $r_{2}(\omega_{i}(1))=\mathrm{E^{RR}}[\omega_{i}(1+m)|\omega_{i}(1)]=\phi^{(1)}(\omega_{i}(t))$,
where we have set $\phi^{(1)}(\omega)=\tau_{0}^{(m-1)}(\tau_{1}^{(1)}(\omega))=\omega\alpha_{m}+\psi_{m}$
with \emph{$\alpha_{m}=\delta_{1}\delta_{0}^{m-1}$} and $\psi_{m}=p_{01}^{(1)}\delta_{0}^{m-1}+p_{01}^{(0)}\frac{1-\delta_{0}^{m-1}}{1-\delta_{0}}$,
and where $\tau_{0}^{(k)}(\omega)=\tau_{0}^{(1)}(\tau_{0}^{(k-1)}(\omega))$
indicates the belief of a node after $k$ slots of passivity when
the initial belief is $\omega$ (i.e., $\tau_{0}^{(k)}(\omega)$ is
obtained recursively by applying $\tau_{0}^{(1)}(\omega)$ to itself
$k$ times). In general, we can obtain $r_{j}(\omega_{i}(1))\mathrm{=E^{RR}}[\omega_{i}(1+(j-1)m)|\omega_{i}(1)]$,
for $j\geq2$, by iterating the procedure above by applying $\phi^{(1)}(\omega)$
to itself $(j-1)$ times. After a little algebra we get $\phi^{(j-1)}(\omega)=\phi^{(1)}(\phi^{(j-2)}(\omega))=\omega\alpha_{m}^{j-1}+\psi_{m}\frac{1-\alpha_{m}^{j-1}}{1-\alpha_{m}}$,
so that $r_{j}(\omega_{i}(1))=\phi^{\left(j-1\right)}(\omega_{i}(1))$,
where we set $\phi^{(0)}(\omega)=\omega$. The reasoning above can
be applied when starting from any arbitrary slot $t$.

Finally, the total reward accrued by the CC from a node that is scheduled
$H$ times, when its belief at the first slot in which it is scheduled
is $\omega$, can be calculated by summing up the average reward $r_{j}(\cdot)$
during each slot in which the node is scheduled (see definition above),
as 
\begin{equation}
\theta^{(H)}\left(\omega\right)=\sum_{j=1}^{H}\beta^{(j-1)m}r_{j}(\omega)=\frac{\psi_{m}}{1-\alpha_{m}}\left(\frac{1-\beta^{mH}}{1-\beta^{m}}-\frac{1-\left(\beta^{m}\alpha_{m}\right)^{H}}{1-\beta^{m}\alpha_{m}}\right)+\frac{1-\left(\beta^{m}\alpha_{m}\right)^{H}}{1-\beta^{m}\alpha_{m}}\omega.\label{eqn: Theta_function}
\end{equation}
Note that, for a node $U_{i}\in\mathcal{G}_{g}$, for $g\geq1$ and
with belief equal to $\omega$ at $t=1$, the first slot in which
the node is scheduled is $t=g$ , and thus its belief at time $t=g$
becomes $\tau_{0}^{\left(g-1\right)}\left(\omega\right)$ (i.e., after
$\left(g-1\right)$ slots of passivity while other groups are scheduled).
Therefore, for a node $U_{i}\in\mathcal{G}_{g}$, with initial belief
$\omega$, the total contribution to the throughput is given by $\beta^{g-1}\theta^{(H)}\left(\tau_{0}^{\left(g-1\right)}\left(\omega\right)\right)$.

Let us now focus on the second step, i.e., proving the inequality
(\ref{eqn: inequality_x_and_y}). At $t=T$, it is easily seen to
hold due to (\ref{eqn: reward_belief_MDP_line_1}) and (\ref{eq: value_DP}).
We then need to show that (\ref{eqn: inequality_x_and_y}) also holds
at $t$. To do so, let us denote as $\mathcal{L}$ and $\mathcal{R}$
the RR policies whose throughputs are given by the LHS and RHS of
(\ref{eqn: inequality_x_and_y}) respectively. The differences between
$\mathcal{L}$ and $\mathcal{R}$ are the positions of the nodes with
belief $x$ and $y$ in the initial belief vectors. Therefore, some
of the $m$ groups created by the two policies might have different
nodes (see the RR operations in Proposition \ref{Prop: MP_Structure}).
To simplify, we refer to the node with belief $x$ ($y$) as node
$x$ ($y$). Let us assume that nodes $x$ and $y$ belong to groups
$\mathcal{G}_{g'}$ and $\mathcal{G}_{g''}$ under policy $\mathcal{R}$,
respectively, while they belong to groups $\mathcal{G}_{g''}$ and
$\mathcal{G}_{g'}$ under policy $\mathcal{L}$, respectively, with
$g''\geq g'$, and $g',\: g''\in\{1,...,m\}$. If $g''=g'$, then
the two policies coincide and (\ref{eqn: inequality_x_and_y}) holds
with equality. If $g''=g'+1$ (nodes are adjacent but do not belong
to the same group), the only difference between policies $\mathcal{L}$
and $\mathcal{R}$ is the scheduling order of nodes $x$ and $y$.

To verify that inequality (\ref{eqn: inequality_x_and_y}) holds,
we need to prove that scheduling node $y$ in group $\mathcal{G}_{g'}$
and node $x$ in group $\mathcal{G}_{g''}$ is no better than doing
the opposite for any $x\geq y$. To elaborate, let $H_{x}^{\mathcal{R}}(t)=H_{y}^{\mathcal{L}}(t)$
and $H_{y}^{\mathcal{R}}(t)$ $=H_{x}^{\mathcal{L}}(t)$ be the number
of times that node $x$ (or $y$) is scheduled under policy $\mathcal{R}$
(or $\mathcal{L}$) and node $y$ (or $x$) is scheduled under policy
$\mathcal{L}$ (or $\mathcal{R}$) in the horizon $\{t,t+1,...,T\}$,
respectively. By recalling (\ref{eqn: Theta_function}) and the discount
factor $\beta$, the contribution generated by node $x$ and $y$
under policy $\mathcal{R}$ is $\beta^{g'-1}\theta^{\left(H_{x}^{\mathcal{R}}(t)\right)}(\tau_{0}^{(g'-1)}(x))$
and $\beta^{g''-1}\theta^{\left(H_{y}^{\mathcal{R}}(t)\right)}(\tau_{0}^{(g''-1)}(y))$
respectively, and similarly under policy $\mathcal{L}$ we have $\beta^{g''-1}\theta^{\left(H_{x}^{\mathcal{L}}(t)\right)}(\tau_{0}^{(g''-1)}(x))$
and $\beta^{g'-1}\theta^{\left(H_{y}^{\mathcal{L}}(t)\right)}(\tau_{0}^{(g'-1)}(y))$.
Note that, in the argument of function $\theta^{(\cdot)}(\cdot)$,
we have considered that the nodes in group $\mathcal{G}_{g'}$ are
scheduled for the first time at slot $g'-1$, and thus the belief
must be updated through function $\tau_{0}^{\left(g'-1\right)}\left(\cdot\right)$,
and similarly for nodes in $\mathcal{G}_{g''}$ the first slot is
$g''-1$. Moreover, the discount factor is $\beta^{g'-1}$ is common
to all the nodes in group $\mathcal{G}_{g'}$, and so is $\beta^{g''-1}$
for group $\mathcal{G}_{g''}$.

By recalling that all the nodes, except $x$ and $y$, are scheduled
at the same slot under the two policies $\mathcal{R}$ and $\mathcal{L}$
(thus giving the same contribution to the throughput), the inequality
(\ref{eqn: inequality_x_and_y}) can thus be reduced to $\beta^{g'-1}\theta^{\left(H_{x}^{\mathcal{R}}(t)\right)}(\tau_{0}^{(g'-1)}(x))+\beta^{g''-1}\theta^{\left(H_{y}^{\mathcal{R}}(t)\right)}(\tau_{0}^{(g''-1)}(y))-\beta^{g''-1}\theta^{\left(H_{x}^{\mathcal{L}}(t)\right)}(\tau_{0}^{(g''-1)}(x))-\beta^{g'-1}\theta^{\left(H_{y}^{\mathcal{L}}(t)\right)}(\tau_{0}^{(g'-1)}(y))\geq0$,
which must hold for all admissible $H_{x}^{\mathcal{R}}(t)=H_{y}^{\mathcal{L}}(t)$
and $H_{y}^{\mathcal{R}}(t)$ $=H_{x}^{\mathcal{L}}(t)$ and all $g''\geq g'$,
with $g',\: g''\in\{1,...,m\}$. There are two cases: \textbf{1})
$H_{x}^{\mathcal{R}}(t)=H_{y}^{\mathcal{L}}(t)=H_{y}^{\mathcal{R}}(t)$
$=H_{x}^{\mathcal{L}}(t)=H\geq1$, that is, nodes $x$ and $y$ are
scheduled the same number of times within the horizon of interest
under the two policies $\mathcal{R}$ and $\mathcal{L}$;\textbf{
2}) $H_{x}^{\mathcal{R}}(t)=H_{y}^{\mathcal{L}}(t)=H$, and $H_{y}^{\mathcal{R}}(t)$
$=H_{x}^{\mathcal{L}}(t)=H-1$, for $H\geq1$, namely, node $x$ (or
$y$) is scheduled one time more than node $y$ (or $x$ ) under policy
$\mathcal{R}$ (or $\mathcal{L}$). By exploiting the RHS of (\ref{eqn: Theta_function}),
after a little algebra, one can verify that the inequality above holds
in both cases, which concludes the proof of Theorem \ref{Thr: Theorem_1}.

\section{Proof of Lemma\label{App: proof_extrema_inequality_m_0_1}\ref{Lem: Extrema_inequalities_lemma}}

\textbf{\textit{\emph{Proof of case a}}}\textbf{)}. From (\ref{eqn: passive_action})-(\ref{eqn: active_action}),
and recalling that $\tau_{0}^{\left(1\right)}\left(0\right)=p_{01}$
from (\ref{eqn: tau_h_whittle}), the leftmost inequality in (\ref{eqn: extrema_inequality_1_m_0_1}.1)
follows immediately as it becomes $V_{m}(0|1)=\beta V_{m}^{\ast}(p_{01})\leq m+\beta V_{m}^{\ast}(p_{01})=V_{m}(0|0)$.
For the rightmost inequality in (\ref{eqn: extrema_inequality_1_m_0_1}.1),
we have $V_{m}(1|1)=1+\beta V_{m}^{\ast}(0)$, while from (\ref{eqn: optimality_equation_SAB})
and the fact that $V_{m}(0|1)\leq V_{m}(0|0)$ we have $V_{m}^{\ast}(0)=\max\left\{ V_{m}(0|0),V_{m}(0|1)\right\} =V_{m}(0|0)$.
Therefore, we have $V_{m}(1|1)=1+\beta V_{m}^{\ast}(0)1+\beta V_{m}(0|0)\geq V_{m}(0|0)$,
which holds as $1+\beta V_{m}(0|0)\geq V_{m}(0|0)$ implies $V_{m}(0|0)\leq\frac{1}{1-\beta}$.
The latter bound always holds, since for $m<1$ the infinite horizon
throughput is upper bounded as $V_{m}^{\ast}(\omega)\leq\sum_{t=0}^{\infty}\beta=\frac{1}{1-\beta}$
given that we can get at most a reward of $R_{m}(\omega,u)\leq1$
in each slot. Hence, inequalities (\ref{eqn: extrema_inequality_1_m_0_1}.1)
are proved. Inequality (\ref{eqn: extrema_inequality_1_m_0_1}.2)
can be proved by contradiction. Specifically, let us assume that:
\emph{hp.1}) $V_{m}(1|0)\geq V_{m}(1|1)$. From (\ref{eqn: optimality_equation_SAB})
we would have $V_{m}^{\ast}(1)=\max\left\{ V_{m}(1|0),V_{m}(1|1)\right\} =V_{m}(1|0)$,
i.e., the passive action would be optimal when $\omega=1$. Moreover,
from (\ref{eqn: passive_action}) we would have $V_{m}(1|0)=m+\beta V_{m}^{\ast}(1)=m+\beta V_{m}(1|0),$
which can be solved with respect to $V_{m}(1|0)$ to get $V_{m}(1|0)=\frac{m}{1-\beta}=V_{m}^{\ast}(1)$.
Therefore, if hypothesis \emph{hp.1}) holds, we also have that $V_{m}(1|1)=1+\beta V_{m}^{\ast}(0)\leq V_{m}(1|0)=V_{m}^{\ast}(1)=\frac{m}{1-\beta}$.
However, the value function $V_{m}^{\ast}(\omega)$ is bounded $\frac{m}{1-\beta}\leq V_{m}^{\ast}(\omega)\leq\frac{1}{1-\beta}$,
where the lower bound is obtained considering a policy that always
chooses the passive action for any belief $\omega$. The boundedness
of the value function, thus implies that if \emph{hp.1}) holds then
$1+\beta\frac{m}{1-\beta}\leq1+\beta V_{m}(0)=V_{m}(1|1)\leq V_{m}(1|0)=\frac{m}{1-\beta}$,
which yields $1+\beta\frac{m}{1-\beta}\leq\frac{m}{1-\beta}$ and
thus $\left(1-\beta\right)\left(1-m\right)\leq0$. But this is clearly
impossible as $m,\beta<1$. Consequently, we have proved that $V_{m}(1|1)\geq V_{m}(1|0)$.

\textbf{\textit{\emph{Proof of case b}}}) Inequality $V_{m}(0|0)\leq V_{m}(0|1)$
follows immediately since $m+\beta V_{m}^{\ast}(p_{01})\leq\beta V_{m}^{\ast}(p_{01})$
holds for $m<0$. The second inequality $V_{m}(0|1)\leq V_{m}(1|1)$
becomes $V_{m}(0|1)\leq=1+\beta V_{m}^{\ast}(0)1+\beta V_{m}(0|1)$,
which leads to $V_{m}(0|1)\leq\frac{1}{1-\beta},$ which always holds
as discussed above. Inequality $V_{m}(1|0)\leq V_{m}(1|1)$ holds
since an active action is always optimal when $m<0$.

\textbf{\textit{\emph{Proof of case c}}}) The inequality holds since
a passive action is always optimal for any $m\geq1$.

\section{Proof of Theorem \ref{thm:indexability}\label{Sec: Proof_indexability}}

Following the discussion in Sec. \ref{eqn: Definition_indexability},
to prove indexability it is sufficient to show that the threshold
$\omega^{\ast}(m)$ is monotonically increasing with the subsidy $m$,
for $0\leq m<1$. In fact, from Proposition \ref{Prop: Threshold_policy}
the passive set (\ref{eqn: Passive_set}) for $m<0$ is $\mathcal{P}(m)=\emptyset$,
while for $m\geq1$, we have $\mathcal{P}(m)=[0,1]$. We then only
need to prove the monotonicity of $\omega^{\ast}(m)$ for $0\leq m<1$,
which has been shown to hold in \cite[Lemma 9]{art:Zhao-whittle}
if
\begin{equation}
\left.\frac{dV_{m}(\omega|1)}{dm}\right\vert _{\omega=\omega^{\ast}(m)}<\left.\frac{dV_{m}(\omega|0)}{dm}\right\vert _{\omega=\omega^{\ast}(m)}.\label{eqn: inequality_monotonicity_lemma}
\end{equation}
 To check if (\ref{eqn: inequality_monotonicity_lemma}) holds, we
differentiate (\ref{eqn: passive_action})-(\ref{eqn: active_action})
at the optimal threshold $\omega=\omega^{\ast}(m)$ as
\begin{eqnarray}
V_{m}(\omega^{\ast}(m)|1) & = & \omega^{\ast}(m)+\beta\omega^{\ast}(m)V_{m}^{\ast}(0)+\beta(1-\omega^{\ast}(m))V_{m}^{\ast}(p_{01}),\;\textrm{and}\label{eqn: conditional_value_fun_at_threshold_activity}\\
V_{m}(\omega^{\ast}(m)|0) & = & m+\beta\left[\tau_{0}^{\left(1\right)}(\omega^{\ast}(m))\left(1+\beta V_{m}^{\ast}(0)\right)+\beta(1-\tau_{0}^{\left(1\right)}(\omega^{\ast}(m)))V_{m}^{\ast}(p_{01})\right],\label{eqn: conditional_value_fun_at_threshold_passivity}
\end{eqnarray}
where (\ref{eqn: conditional_value_fun_at_threshold_passivity}) follows
from (\ref{eqn: active_action}) and from the fact that $\tau_{0}^{(1)}(\omega)\geq\omega$,
for any $\omega$ (see (\ref{eqn: tau_h_whittle})), and hence $V_{m}^{\ast}(\tau_{0}^{(1)}(\omega^{\ast}(m)))=V_{m}(\tau_{0}^{(1)}(\omega^{\ast}(m))|1)$,
since arm activation is optimal for any $\omega>\omega^{\ast}(m)$.

By letting $D_{m}(\omega)=\frac{dV_{m}^{\ast}(\omega)}{dm}$, then
from (\ref{eqn: conditional_value_fun_at_threshold_activity}) we
have $\left.\frac{dV_{m}(\omega|1)}{dm}\right\vert _{\omega=\omega^{\ast}(m)}=\beta\omega^{\ast}(m)D_{m}(0)+\beta(1-\omega^{\ast}(m))D_{m}(p_{01})$,
while from (\ref{eqn: conditional_value_fun_at_threshold_passivity})
we get $\left.\frac{dV_{m}(\omega|0)}{dm}\right\vert _{\omega=\omega^{\ast}(m)}=1+\beta^{2}\tau_{0}^{\left(1\right)}(\omega^{\ast})D_{m}(0)+\beta^{2}(1-\tau_{0}^{\left(1\right)}(\omega^{\ast}))D_{m}(p_{01})$.
Finally, after some algebraic manipulations, and recalling that $D_{m}(0)=\frac{dV_{m}^{\ast}(0)}{dm}=\frac{d\left(m+\beta V^{\ast}(p_{01})\right)}{dm}=1+recursively\beta D_{m}(p_{01})$,
we can rewrite (\ref{eqn: inequality_monotonicity_lemma}) as\newline
$D_{m}(p_{01})\beta\left(1-\beta\right)\left[1-\mathcal{\omega}\left(1-\beta(1-p_{01})\right)\right]+\beta\left[\mathcal{\omega}\left(1-\beta(1-p_{01})\right)-\beta p_{01}\right]<1$.
To show that the last inequality holds when $0\leq m<1$, we first
upper bound the derivative of the value function as $D_{m}(\omega)$
$\leq\frac{1}{1-\beta}$, since $\frac{d}{dm}R_{m}(\omega)\leq1$.
Finally, using this upper bound $D_{m}(p_{01})$ $\leq\frac{1}{1-\beta}$
after a little algebra (\ref{eqn: inequality_monotonicity_lemma})
reduces to $\beta(1-\beta p_{01})<1,$ which clearly holds for any
$\beta\in[0,1)$ as $0\leq p_{01}\leq1$. This concludes the proof
of Theorem \ref{thm:indexability}.


\begin{thebibliography}{1}

\bibitem{art:Data_networks} D. Bertsekas, R. G. Gallager, \emph{Data
Networks}. Englewood Cliffs, NJ: Prentice Hall, 1992.

\bibitem{art:Grid_computing_2} A. Benoit, L. Marchal, J.-F. Pineau,
Y. Robert, F. Vivien, {}``Scheduling concurrent bag-of-tasks applications
on heterogeneous platforms,'' \emph{IEEE Trans. Computers}, vol.
59, no. 2, pp. 202-217, Feb. 2010.

\bibitem{art:Virtual_Machine}Y. Bai, C. Xu and Z. Li, \textquotedblleft{}Task-aware
based co-scheduling for virtual machine system,\textquotedblright{}
in\emph{ Proc. ACM Symp. On Applied Comp}., Sierre, Switzerland, pp.
181-188, Mar. 2010.

\bibitem{art:Monahan} G. E. Monahan, {}``A survey of partially observable
Markov decision processes: Theory, models, and algorithms,\textquotedblright{}\textit{
Manag. Sci}., vol. 28, no. 1, pp. 1-16, 1982.

\bibitem{book:Gittins} J. Gittins, K, Glazerbrook, R. Weber, \textit{Multi-armed
Bandit Allocation Indices}. West Sussex, UK: Wiley, 2011.

\bibitem{art:whittle}P. Whittle, \textquotedblleft{}Restless bandits:
Activity allocation in a changing world,\textquotedblright{} \emph{J.
Appl. Probab}., vol. 25, pp. 287-298, 1988.

\bibitem{art:Zhao-myopic} S. H. A. Ahmad, M. Liu, T. Javidi, Q. Zhao
and B. Krishnamachari, {}``Optimality of myopic sensing in multi-channel
opportunistic access,\textquotedblright{}\textit{ IEEE Trans. Inf.
Theory}, vol. 55, No. 9, pp. 4040-4050, Sept. 2009.

\bibitem{art:Ahmad-vector-myopic} S. H. A. Ahmad, M. Liu, {}``Multi-channel
opportunistic access: A case of restless bandits with multiple plays,\textquotedblright{}
in \textit{Proc. 47th Ann. Allerton Conf. Commun., Contr., Comput.},
Monticello, IL, pp. 1361-1368, Sept. 2009.

\bibitem{art:Zhao-whittle} K. Liu and Q. Zhao, {}``Indexability
of restless bandit problems and optimality of Whittle index for dynamic
multichannel access,\textquotedblright{}\textit{ IEEE Trans. Inf.
Theory}, vol. 56, no. 11, pp. 5547-5567, Nov. 2010.

\bibitem{art:Cassandra} L. P. Kaelbling, M. L. Littman, and A. R.
Cassandra, {}``Planning and acting in partially observable stochastic
domains,\textquotedblright{}\textit{ Artif. Intell.}, vol. 101, pp.
99-134, May 1998. 

\bibitem{book:Puterman} M. L. Puterman, \textit{Markov Decision Processes:
Discrete Stochastic Dynamic Programming}. Hoboken, NJ: Wiley, 2005.

\bibitem{art:Mora}D. Bertsimas and J. E. Niño-Mora, \textquotedblleft{}Restless
bandits, linear programming relaxations, and a primal-dual heuristic,\textquotedblright{}
\emph{Oper. Res.}, vol. 48, no. 1, pp. 80-90, Jan. 2000.

\end{thebibliography}
\end{document}